\documentclass[12pt]{article}
\usepackage{amssymb,amsmath,amsthm,tikz-cd,amscd,cite,hyperref, pgfplots, authblk,tikz, caption, graphicx}
\usepackage[letterpaper, margin=1in]{geometry}

\title{Heat Kernel Estimates for Schr\"odinger Operators with Decay at Infinity on Parabolic Manifolds}

\author[1]{Anthony Graves-McCleary} 
\author[2]{Laurent Saloff-Coste}

\affil[1]{Department of Mathematics, Cornell University, Ithaca, NY, USA ag2537@cornell.edu}
\affil[2]{Department of Mathematics, Cornell University, Ithaca, NY, USA lsc@math.cornell.edu}
\date{}
\newcommand{\set}[1]{\left\{#1\right\}}
\newcommand{\Real}{\mathbf R}

\newcommand{\N}{\mathbf{N}}

\newcommand{\ip}[2]{\left< #1, #2  \right>}
\newcommand{\brac}[1]{\left< #1 \right>}
\newcommand{\abs}[1]{\left\vert#1\right\vert}
\newcommand{\norm}[1]{\left\Vert#1\right\Vert}

\newcommand{\fct}[3]{#1 \colon #2 \rightarrow #3}

\renewcommand{\phi}{\varphi}

\newcommand{\mc}{\mathcal}

\newcommand{\grad}{\nabla}

\newtheorem{theorem}{Theorem}[section]
\newtheorem{proposition}[theorem]{Proposition}
\newtheorem{lemma}[theorem]{Lemma}
\newtheorem{definition}[theorem]{Definition}
\newtheorem{corollary}[theorem]{Corollary}
\newtheorem{example}[theorem]{Example}

\usetikzlibrary{arrows, matrix}

\numberwithin{equation}{section}

\begin{document}
\maketitle

\begin{abstract}
We give estimates for positive solutions for the Schr\"odinger equation $(\Delta_\mu+W)u=0$ on a wide class of parabolic weighted manifolds $(M, d\mu)$ when $W$ decays to zero at infinity faster than quadratically. These can be combined with results of Grigor'yan \cite{grigoryan2} to give matching upper and lower bounds for the heat kernel of the corresponding Schr\"odinger operator $\Delta_\mu+W$. In particular, this appears to complement known results for Schr\"odinger operators on $\Real^2$.
\end{abstract}

\section{Introduction}\label{introsection}

We are interested in estimating the heat kernel for operators of the form $\Delta_\mu+W$, where $\Delta_\mu$ is the Laplacian of a weighted manifold $(M, d\mu)$ and $W$ is a function that decays faster than quadratically at infinity. We state our results on parabolic manifolds satisfying the parabolic Harnack inequality, a connectivity of annuli condition and a mild condition on volume growth. Notably our main result Theorem \ref{maithm} appears to give an addition to known results on $\Real^2$; see Corollary \ref{plane}.


\subsection{Notation}

For $x=(x_1, \dots, x_n)\in \Real^n$ we let $\abs{x}=\left(\sum_{j=1}^n x_j^2\right)^{1/2}$ denote the norm of $x$. For a Riemannian manifold $M$ we let $d(x,y)$ denote the geodesic distance between $x,y\in M$. When $M$ has a distinguished point $o\in M$, we let $\abs{x}=d(x,o)$ and $\brac{x}=2+\abs{x}$.\\

We let $a\preceq b$ denote the existence of a constant $C>0$ such that $a\leq C\cdot b$. We let $a\asymp b$ denote that both $a\preceq b$ and $a\succeq b$ hold. Exact constants may sometimes change line to line. Given a number $a\in \Real$ we let $a_+=\max(a, 0)$ and $a_-=\max(-a, 0)$. Given a real valued function $f$ we let $f_+=\max(f, 0)$ and $f_-=\max(-f, 0)$ so that $f_+, f_-\geq 0$ and $f=f_+-f_-$. Given a smooth manifold $M$ we let $$C^\infty(M)=\set{\fct{f}{M}{\Real} \textrm{ smooth}}$$ $$C_c^\infty(M)=\set{f\in C^\infty(M)\colon f\textrm{ has compact support}}$$ $$C_0^\infty(M)=\set{f\in C^\infty(M)\colon f\textrm{ vanishes at infinity}}.$$ 
Given a Borel measure $d\mu$ on a Riemannian manifold $M$, we let $B(x,r)$ and $\overline{B}(x,r)$ denote the open and closed geodesic balls of center $x\in M$ and radius $r>0$ respectively. We let $V_\mu(x,r)=\mu(B(x,r))$. If $M$ has a distinguished point $o\in M$, then we let $V_\mu(r)$ denote $V_\mu(o, r)$.

\subsection{Background and Overview}

The problem of estimating the heat kernel of $\Delta +W$ for various functions $W$ has attracted considerable research interest. Aizenman and Simon in \cite{aizenmann1} studied Harnack inequalities for Schr\"odinger operators and proved some heat kernel estimates in well-behaved bounded sets with Dirichlet boundary conditions. Due to the useful concept of Green-boundedness, which we will discuss below in Definition \ref{greenb}, many heat kernel estimates for Schr\"odinger operators have been obtained on manifolds that are \textit{non-parabolic}, see Definition \ref{nonpara}. See for example the papers by Takeda \cite{takeda1}, Devyver \cite{devyver1}, and Chen and Wang \cite{chenwang1} for results on heat kernels of Schr\"odinger operators on non-parabolic manifolds. See Davies and Simon \cite{daviessimon1} and Zhang \cite{zhang2}, \cite{zhang1} for results on Schr\"odinger operators with potentials that decay exactly quadratically. Fraas, Krej\v{c}i\v{r}\'ik, and Pinchover in \cite{pinchover2} studied asymptotics of ratios of heat kernels for both subcritical and critical Schr\"odinger operators. Grigor'yan in \cite{grigoryan2} obtained results for the heat kernel of a Schr\"odinger operator whose potential is positive and compactly supported on a parabolic manifold. Shen in \cite{shen1} studied the Green's functions of Schr\"odinger operators in $\Real^n$, $n\geq 3$. Murata in \cite{murata1} and \cite{murata3} obtained large-time asymptotics for heat kernels of subcritical Schr\"odinger operators on Euclidean spaces including $\Real^2$.\\



Of particular relevance is the following result of Murata. See \ref{definitionprofile} for the definition of the profile and \ref{critdef} for the definition of subcritical and critical. Here we let $\brac{x}=2+\abs{x}$.

\begin{theorem}\label{muratatheorem} (Murata \cite{murata2}) Let $W$ be such that $\brac{x}^{\alpha-\frac{2}{p}}\abs{W(x)}\in L^p(\Real^2, dx)$ for $p>1$, $\alpha>2$. Then if $\Delta+W$ is subcritical, the profile $h>0$ for $\Delta+W$ satisfies $$h(x)\asymp \log \brac{x}.$$ If $\Delta+W$ is critical, the profile $h>0$ for $\Delta+W$ satisfies $$h(x)\asymp 1.$$
\end{theorem}

See also Murata \cite{murata4}. The result below is a product of combining Murata's theorem above with Theorem 10.10 in Grigor'yan \cite{grigoryan2}.

\begin{theorem}\label{muratagrigoryan}
Let $W\in C^\infty(\Real^2)$ be such that $\brac{x}^{\alpha-\frac{2}{p}}\abs{W(x)}\in L^p(\Real^2, dx)$ for $p>1$, $\alpha>2$. Let $p^W(t,x,y)$ denote the heat kernel of $\Delta+W$. There exist constants $c_1, c_2, c_3, c_4>0$ such that:

\noindent(i) If $\Delta+W$ is subcritical, then for all $t>0$, $x,y\in \Real^2$, $$\frac{c_1\log\brac{x}\log \brac{y}}{t\log(\brac{x}+\sqrt{t})\log(\brac{y}+\sqrt{t})}e^{-c_2\frac{\abs{x-y}^2}{t}}\leq p^W(t,x,y)\leq \frac{c_3 \log\brac{x}\log \brac{y}}{t\log(\brac{x}+\sqrt{t})\log(\brac{y}+\sqrt{t})}e^{-c_4\frac{\abs{x-y}^2}{t}}.$$

\noindent(ii) If $\Delta+W$ is critical, then for all $t>0$, $x,y\in \Real^2$, $$\frac{c_1}{t}e^{-c_2\frac{\abs{x-y}^2}{t}}\leq p^W(t,x,y)\leq \frac{c_3}{t}e^{-c_4\frac{\abs{x-y}^2}{t}}.$$
\end{theorem}

In this paper we study Schr\"odinger operators whose potentials are smooth and vanish at infinity, and in that context give a result analogous to Murata's for a wide class of parabolic weighted manifolds. This can be then combined with Grigor'yan's theorem to give heat kernel estimates, although provide our own proof of this in Theorem \ref{maithm}. Looking back to $\Real^2$, our main result offers an addition to Theorem \ref{muratagrigoryan} by including those $W\in C_0^\infty(\Real^2)$ for which $\log\brac{x}\abs{W(x)}\in L^1(\Real^2, dx)$.










\subsection{Setting: Weighted Manifolds}
\begin{definition}
A \textbf{weighted manifold} (also called a manifold with density) is a pair $(M, d\mu)$ where $M$ is a Riemannian manifold without boundary and $d\mu$ is a Borel measure on $M$ that has a smooth positive density with respect to canonical Riemannian volume measure on $M$. 
\end{definition}

Given a smooth function $f$ on $M$, let $\grad f$ denote the gradient of $f$ on $M$. Note that $\grad f$ is defined using only $f$ and the Riemannian metric, and is therefore independent of any choice of weighted measure on $M$. We define $\mc{F}_\mu$ to be the completion of $C_c^\infty(M)$ under the norm $$\norm{f}_{\mc{E}^1}=\int_M \abs{f}^2d\mu+\int_M \abs{\grad f}^2d\mu.$$ There exists a (minimal) operator $\Delta_\mu$ called the \textbf{Laplacian} of $(M, d\mu)$ such that whenever $f,g\in C_c^\infty(M)$, $$\int_M \ip{\grad f}{\grad g}d\mu = \int_M f\Delta_\mu (g) d\mu.$$ Here $\ip{\cdot}{\cdot}$ is the Riemannian inner product. \\

\noindent\textbf{Remark:} We take the convention that $\Delta_\mu\geq 0$, so that e.g. on $\Real^n$ with the Euclidean metric and $d\mu=$ Lebesgue measure, we have $\Delta_\mu=-\sum_{j=1}^n\frac{\partial^2}{\partial x_j^2}$.\\

 Continuing, the quadratic form \begin{equation}\label{dirichletenergy}\mc{E}_\mu(f,f)=\int_M \abs{\grad f}^2d\mu\end{equation} is a strictly local Dirichlet form on $L^2(M, d\mu)$ with domain $\mc{F}_\mu$. This Dirichlet form has infinitesimal generator $\Delta_\mu$. We have the associated heat semigroup $e^{-t\Delta_\mu}$ and ensuing heat kernel $p_\mu(t,x,y)$, which is a fundamental solution of the heat equation $(\partial_t+\Delta_\mu)u=0$.

\begin{definition}\label{nonpara}
Let $(M, d\mu)$ be a complete, non-compact weighted manifold and let $p_\mu(t,x,y)$ be its heat kernel. If there exist $x,y\in M$ distinct such that $$\int_0^\infty p_\mu(t,x,y)dt=+\infty,$$ we say $(M, d\mu)$ is \textbf{parabolic}. Otherwise we say $(M, d\mu)$ is \textbf{non-parabolic}.
\end{definition}

For example, $(\Real^n,dx)$ is parabolic for $n=1, 2$ and non-parabolic for $n\geq 3$.

\begin{definition}
Let $(M, d\mu)$ be a complete, non-compact, non-parabolic weighted manifold. The \textbf{Green's function} $G_\mu$ of $(M, d\mu)$ is defined as follows: for $x,y\in M$ distinct, $$G_\mu(x,y)=\int_0^\infty p_\mu(t,x,y)dt.$$
\end{definition}

The Green's function is the fundamental solution of $\Delta_\mu$, and for each $x\in M$, the function $G_\mu(x, \cdot)$ is $\Delta_\mu$-harmonic on $M\setminus \set{x}$, which we define below.

\begin{definition}
Let $\Omega\subseteq M$ be open and let $u$ be a smooth function on $\Omega$. We say that $u$ is \textbf{$\Delta_\mu$-harmonic} on $\Omega$ if $$\Delta_\mu u=0.$$ We say that $u$ is \textbf{$\Delta_\mu$-subharmonic} (resp. \textbf{$\Delta_\mu$-superharmonic}) on $\Omega$ if $$\Delta_\mu u\leq 0 \hspace{5mm} \mathrm{(resp. }\textrm{ } \Delta_\mu u\geq 0\mathrm{)}.$$ 
\end{definition}

\begin{definition}\label{definitionprofile}
Let $(M, d\mu)$ be a complete, non-compact weighted manifold and let $W\in C^\infty(M)$. A \textbf{profile} for $\Delta_\mu+W$ is a positive smooth function $h>0$ satisfying $$(\Delta_\mu+W)h=0.$$

\end{definition}

\begin{definition}
Let $(M, d\mu)$ be a complete, non-compact weighted manifold and let $W\in C^\infty(M)$. We say that $\Delta_\mu+W\geq 0$ if for all $\phi\in C_c^\infty(M)$ we have $$\int_M \phi(\Delta_\mu+W)\phi d\mu\geq 0.$$
\end{definition}

We finish this section by citing two useful classical results.

\begin{theorem}\label{allegretto} (Allegretto-Piepenbrink Theorem) Let $(M, d\mu)$ be a complete, non-compact weighted manifold and let $W\in C^\infty(M)$ be bounded. Then $\Delta_\mu+W\geq 0$ if and only if there exists a profile for $\Delta_\mu+W$.

\end{theorem}

\noindent\textbf{Proof:} See e.g. the book by Pigola, Rigoli, and Setti \cite{pigola1}, Lemma 3.10 for unweighted manifolds. The proof for weighted manifolds is the same.\qed

\begin{theorem}\label{liouville}
(Liouville's Theorem) Let $(M, d\mu)$ be a complete, non-compact, parabolic weighted manifold. Let $u>0$ be a positive $\Delta_\mu$-superharmonic function on $M$. Then $u$ is constant.
\end{theorem}

\noindent\textbf{Proof:} See e.g. Grigor'yan \cite{grigoryan3} Theorem 5.1 for the case when $M$ is unweighted, or Sturm \cite{sturm1} Theorem 3 for a more general version in a strictly local recurrent Dirichlet space.\qed

\section{Parabolic Manifolds and Harmonic Functions Outside a Ball}\label{outsideball}

Let $(M, d\mu)$ be a complete, non-compact, parabolic weighted manifold. Let $d$ denote geodesic distance on $M$. We let $B(x,r)$ denote the ball of radius $r>0$ and center $x\in M$, and furthermore let $V_\mu(x,r)=\mu(B(x,r))$. When $M$ has a distinguished point $o$, we let $V_\mu(r)=V_\mu(o, r)$.

\begin{definition}
We say that a Riemannian manifold $M$ satisfies the \textbf{relatively connected annuli (RCA) property with respect to $o\in M$} if there exists $A>1$ such that for all $x,y\in M$ and all $R>0$, if $d(x,o)=d(y,o)=R$, there exists a continuous path from $x$ to $y$ in the annulus $B(o, AR)\setminus B(o, A^{-1}R)$.
\end{definition} 

For example, $\Real^2$ satisfies (RCA) with respect to the origin (or any point), and $\Real$ fails to satisfy (RCA). A connected sum of two Euclidean spaces fails to satisfy (RCA) with respect to any point.\\

\begin{definition}
We say that a complete, non-compact weighted manifold $(M, d\mu)$ satisfies \textbf{Gaussian-type heat kernel estimates (HKE)} if there exist constants $c_1, c_2, c_3, c_4>0$ such that for all $t>0$ and all $x,y\in M$ we have \begin{equation}
\frac{c_1}{V_\mu(x,\sqrt{t})}\exp\left(-c_2\frac{d(x,y)^2}{t}\right)\leq p_\mu(t,x,y)\leq \frac{c_3}{V_\mu(x,\sqrt{t})}\exp\left(-c_4\frac{d(x,y)^2}{t}\right).
\end{equation} Here $p_\mu(t,x,y)$ is the heat kernel of $(M, d\mu)$ and $V_\mu(x,r)=\mu(B(x,r))$.
\end{definition}

The property (HKE) has wide-ranging implications, some of which we will discuss in \ref{laurentbigthm}. An application to Green's functions that will be particularly useful to us is stated below.

\begin{lemma}\label{greenhke}
Let $(M, d\mu)$ be a complete, non-compact weighted manifold satisfying (HKE). Then the Green's function $G_\mu(x,y)$ of $(M, d\mu)$ satisfies \begin{equation}
G_\mu(x,y)\asymp \int_{d(x,y)^2}^\infty \frac{dt}{V_\mu(x,\sqrt{t})}.
\end{equation}
\end{lemma}

\noindent\textbf{Proof:} See e.g. Saloff-Coste \cite{aspects1}, Corollary 5.4.13.\qed\\

Now consider the following two properties.

\begin{definition}
A weighted manifold $(M, d\mu)$ satisfies \textbf{volume doubling (VD)} if there exists a constant $C_{VD}>0$ such that for all $x\in M$ and $r>0$, $\mu(B(x, 2r))\leq C\mu(B(x,r))$.
\end{definition}

\begin{definition}
A weighted manifold $(M, d\mu)$ satisfies the \textbf{Poincar\'e inequality (PI)} if there exist constants $C_P, A_P>0$ such that for all $x\in M$, $r>0$, and for all $f\in C^1(B(x, A_Pr))$, we have \begin{equation}
\int_{B(x,r)}\abs{f-f_{B(x,r)}}^2d\mu\leq C_Pr^2\int_{B(x, Ar)}\abs{\grad f}^2d\mu.
\end{equation} Here $f_B=\frac{1}{\mu(B)}\int_B fd\mu$.
\end{definition}


We now state an important theorem connecting (HKE) to (VD) and (PI).

\begin{theorem} \label{laurentbigthm}(Grigor'yan \cite{grigoryan1}, Saloff-Coste \cite{saloffcoste1}) Let $(M, d\mu)$ be a complete, non-compact weighted Riemannian manifold. Then $(M, d\mu)$ satisfies (HKE) if and only if $(M, d\mu)$ satisfies (VD) and (PI).

\end{theorem}

In practice, (VD) and (PI) are often easier to check than (HKE). Gaussian-type heat kernel estimates imply a wide array of properties, including the parabolic Harnack inequality. The parabolic Harnack inequality, which is actually equivalent to (HKE), then implies the (scale-invariant) elliptic Harnack inequality, which we state now.

\begin{lemma}\label{ehi}
Let $(M, d\mu)$ be a complete, non-compact weighted manifold satisfying (HKE). There exists a constant $C_H>0$ such that for any open ball $B=B(x,r)$ in $M$ and any positive $\Delta_\mu$-harmonic function $u>0$ on $2B=B(x, 2r)$, we have $$\sup_B u \leq C_H \inf_B u.$$ 
\end{lemma}

\noindent\textbf{Proof:} See e.g. Saloff-Coste \cite{saloffcoste1} Theorem 3.1.\qed\\

 The conclusion of the following classical lemma is known as the metric doubling property.

\begin{lemma}\label{metricdoubling}
Assume that $(M, d\mu)$ is a complete, non-compact, weighted manifold satisfying (VD). There exists $N\in \N$ such that for all $x\in M$ and all $R>0$, the ball $B(x, R)$ of radius $R$ and center $x$ can be covered by at most $N$ balls of radius at most $R/2$.

\end{lemma}

\noindent\textbf{Proof:} See e.g. the book by Coifman and Weiss \cite{coifmanweiss}, pages 67-68.\qed\\

The above lemma can be combined with the scale-invariant elliptic Harnack inequality to produce the following radial version of the Harnack inequality.

\begin{lemma}\label{radialh} Assume that $(M, d\mu)$ is a complete, non-compact, weighted manifold satisfying (HKE) and (RCA) with respect to $o\in M$. Let $A>0$ be the constant appearing the definition of (RCA). Let $R>0$ and let $u>0$ be $\Delta_\mu$-harmonic on $B(o, 2AR)\setminus \overline{B}(o, (2A)^{-1}R)$. Then there exists $C>0$ independent of $R$ such that if $d(x, o)=d(y, o)$, then $u(x)\leq Cu(y)$. 

\end{lemma}

\noindent\textbf{Proof:} By iterating Lemma \ref{metricdoubling}, we get $N\in \N$ such that for any $R>0$ and any $x\in M$, the ball $B(x, R)$ can be covered in $N$ balls of radius at most $R/16A^2$. Now fix $R>0$ and cover the ball of radius $B(o, 2AR)$ with $N$ balls of radius $R/8A$. Note that if $B$ is such a ball that meets $B(o, AR)\setminus \overline{B}(o, A^{-1}R)$, then $2B$ lies in $B(o, 2AR)\setminus \overline{B}(o, (2A)^{-1}R)$.\\

Let $x$ and $y$ be such that $d(x, o)=d(y, o)=R$. By (RCA), there exists a continuous path from $x$ to $y$ in $B(o, A)\setminus \overline{B}(o, A^{-1}R)$. By the previous, our path is covered by at most $N$ balls $B$ such that $2B\subseteq B(o, 2A)\setminus\overline{B}(o, (2A)^{-1}R)$. By Lemma \ref{ehi}, there exists a constant $C_H>0$ independent of $R$ such that for each such ball $B$ we have \begin{equation}\label{ellipticharnack}\sup_B u\leq C_H\inf_B u.\end{equation} Therefore by applying the \ref{ellipticharnack} to each of the balls covering the path without repetition, we obtain the desired result, namely that $u(x)\leq (C_H)^N u(y)$.\qed\\

Let $(M, d\mu)$ be a weighted manifold with distinguished point $o\in M$. For $r>0$ let $V_\mu(r)=V_\mu(o, r)$ be the volume of the ball $B(o, r)$. Given $x\in M$ we set $\abs{x}:=d(x,o)$.\\

Let us define \begin{equation}
H(r):=1+\left(\int_1^{r^2} \frac{ds}{V_\mu(\sqrt{s})}\right)_+.
\end{equation}

Note that $$H(r)\asymp 1+\left(\int_1^{r} \frac{s}{V_\mu(s)}ds\right)_+\asymp\int_0^{r} \frac{se^{-1/s}}{V_\mu(s)}ds.$$

We wish to study positive functions that are harmonic outside a neighborhood of $o$. We have the following result due to Grigory'an and Saloff-Coste \cite{grigoryansaloffcoste1}:

\begin{proposition} \label{dirharmonic}
Let $(M, d\mu)$ be a complete, non-compact, parabolic weighted manifold satisfying (HKE) as well as (RCA) with respect to $o\in M$. Let $K$ be compact with nonempty interior and $o\in K$. There exists a positive smooth function $u$ on $M$ which is $\Delta_\mu$-harmonic in $M\setminus K$ and satisfies \begin{equation}
u(x)\asymp H(\abs{x})
\end{equation} for all $x\in M$.
\end{proposition}

Note that because $M$ is parabolic, $H(\abs{x})\rightarrow +\infty$ as $\abs{x}\rightarrow \infty$.\\

We now introduce a potential $q$, which we will require to be well-behaved: namely smooth, nonnegative, and compactly supported. Profiles, defined in \ref{definitionprofile} as positive solutions $h>0$ to $(\Delta_\mu+q)h=0$ will be highly important to our main results.\\

\begin{lemma}
Let $(M, d\mu)$ be a complete, non-compact, parabolic weighted manifold satisfying (HKE) as well as (RCA) with respect to $o\in M$. Let $q\in C_c^\infty(M)$ be non-negative and not identically zero. Let $h>0$ satisfy $(\Delta_\mu+q)h=0$. Then $$h(x)\asymp H(\abs{x}).$$
\end{lemma}

\noindent\textbf{Proof:} Let $A>0$ be as in the definition of (RCA) and let $h>0$ be a profile for $\Delta_\mu+q$ as in the hypothesis. For $x\in M$ let $\abs{x}=d(x,o)$. Let $R_0>0$ be such that $B(o, (2A)^{-1}R_0)$ contains the support of $q$. Let $K=\overline{B}(o, (2A)^{-1}R_0)$. Note that $h$ is $\Delta_\mu$-harmonic in $M\setminus K$.\\

As in Proposition \ref{dirharmonic}, let $u>0$ be smooth on $M$ that is $\Delta_\mu$-harmonic in $M\setminus K$ such that $u(x)\asymp H(\abs{x})$ for all $x\in M$. By Lemma \ref{radialh}, there exists $C>0$ such that if $\abs{x}=\abs{y}>R_0$, then $u(x)\leq Cu(y)$, and similarly for $h$. With this in mind for $r>0$ and a function $f$ on $M$ we define $f(r)=\sup_{\abs{x}=r}f(x)$.\\

It suffices to show that $h(x)\asymp u(x)$. Since $(\Delta_\mu+q)h=0$, $h$ is non-constant. Moreover we have $\Delta_\mu h=-qh\leq 0$ and so $h$ is $\Delta_\mu$-subharmonic. If $h$ were bounded, then for some $D>0$, $D-h$ would be a positive $\Delta_\mu$-superharmonic function in $M$ which contradicts Liouville's Theorem \ref{liouville}. So $h$ is unbounded and in particular $\limsup_{r\rightarrow \infty} h(r)=+\infty$.\\

We can thus choose $\epsilon>0$ sufficiently small so that for some $R>R_0$, $(h-\epsilon u)(R)>(h-\epsilon u)(R_0)>0$. By the maximum principle, we get that $(h-\epsilon u)(r)\geq (h-\epsilon u)(R)$ for all $r>R$. Thus if $\abs{x}>R$, then let $y$ be such that $\abs{y}=\abs{x}$ and $(h-\epsilon u)(y)=(h-\epsilon u)(r)$. By Harnack's inequality \ref{radialh} we have that $u(x)\leq Cu(y)\leq \frac{C}{\epsilon} h(y)\leq \frac{C^2}{\epsilon}h(x)$. Thus we have $u\preceq h$ outside a compact set. This implies $u\preceq h$ on all of $M$.\\

Since all we needed for this argument was $h$ and $u$ to be positive, unbounded, and $\Delta_\mu$-harmonic in $M\setminus K$, we can repeat this argument with $h$ and $u$ switched, yielding $h\asymp u$.\qed

\section{Doob's $h$-transform}\label{htransformsection}

Let $(M, d\mu)$ be a weighted manifold with Laplacian $\Delta_\mu$. Let $W\in C^\infty(M)$ be bounded. Consider the Schr\"odinger form on $M$ given by \begin{equation}
\mc{E}^W_\mu(f,f)=\int_M \abs{\grad f}^2d\mu+\int_M W\abs{f}^2d\mu
\end{equation} with domain $\mc{F}_\mu$, the same as that of $\mc{E}_\mu$ from \ref{dirichletenergy}. Assume there exists a profile $h>0$ for $\Delta_\mu +W$ and let $d\nu=h^2d\mu$. Consider the weighted manifold $(M, d\nu)$ with strictly local Dirichlet form \begin{equation}
\mc{E}_\nu(f,f)=\int_M \abs{\grad f}^2d\nu.
\end{equation}

We can contexualize this in terms of the unitary operator $\fct{\mathfrak{u}_h}{L^2(M, h^2d\mu)}{L^2(M, d\mu)}$ given by $\mathfrak{u}_h(f)=hf$. Integration by parts and the fact that $h$ is a profile for $\Delta_\mu+W$ imply that \begin{equation}\label{intbypartsdoob}\mc{E}_\nu(f, f)=\mc{E}_\mu^W(\mathfrak{u}_h(f), \mathfrak{u}_h(f))\end{equation} for $f\in C_c^\infty(M)$. This suggests defining the domain of $\mc{E}_\nu$ as $\mc{F}_\nu=\mathfrak{u}_h^{-1}\mc{F}_\mu$.\\

Let $d\nu=h^2d\mu$. By Equation \ref{intbypartsdoob}, the Laplacian of the weighted manifold $(M, d\nu)$ is related to that of $(M, d\mu)$ by $\Delta_\nu(f)=\frac{1}{h}(\Delta_\mu+W)(hf)$, i.e. $\Delta_\nu=\mathfrak{u}_h^{-1}\circ (\Delta_\mu+W)\circ \mathfrak{u}_h$. Let $p_\nu(t,x,y)$ denote the heat kernel of $\Delta_\nu$. One can verify that $p_\nu(t,x,y)$ is related to the heat kernel $p^W_\mu(t,x,y)$ of $\Delta_\mu+W$ via \begin{equation}\label{doobhke}
p^W_\mu(t,x,y)=h(x)h(y)p_\nu(t,x,y).
\end{equation} Thus, an estimate of $p_\nu(t,x,y)$ together with an estimate of $h$ yield an estimate of $p^W_\mu(t,x,y)$.

\subsection{Iterated $h$-Transform}

In this section we discuss the iterated $h$-transform. 

\begin{lemma}\label{iter} Let $(M, d\mu)$ be a weighted manifold and let $W_1, W_2\in C^\infty(M)$. Let $g,h>0$ be positive functions on $M$ satisfying $(\Delta_\mu+W_1)h=0$ and $(\Delta_\nu+W_2)g=0$, where $d\nu=h^2d\mu$. Then letting $d\theta=g^2h^2d\mu$, if $p^{W_1+W_2}_\mu(t,x,y)$ denotes the heat kernel of $\Delta_\mu+W_1+W_2$ on $(M, d\mu)$ and $p_\theta(t,x,y)$ denotes the heat kernel of $\Delta_\theta$ on $(M, d\theta)$, we have that \begin{equation}
p^{W_1+W_2}_\mu(t,x,y)=g(x)h(x)g(y)h(y)p_\theta(t,x,y)
\end{equation} for all $t>0$ and $x,y\in M$.
\end{lemma}

\noindent\textbf{Proof:} Note that \begin{equation}
(\Delta_\mu+W_1+W_2)(gh)=h\Delta_\nu g+W_2gh=h(\Delta_\nu +W_2)g=0.
\end{equation} Thus $gh$ is a profile for $\Delta_\mu+W_1+W_2$. The conclusion then follows from Equation \ref{doobhke}.\qed\\

Thus an estimate on $p_\theta(t,x,y)$ yields an estimate on $p^{W_1+W_2}_\mu(t,x,y)$. In practice, we will choose $W_1$ so that the profile $h$ has tractible form and such that $(M, h^2d\mu)$ satisfies (HKE). We will then choose $W_2$ so that $g\asymp 1$ and such that $W_1+W_2$ is of the desired form, leading to a heat kernel estimate for $\Delta_\mu+W_1+W_2$.\\

\begin{lemma}\label{itergreen}
Let $(M, d\mu)$ be a weighted manifold and let $W_1, W_2\in C^\infty(M)$. Let $h>0$ be smooth satsifying $(\Delta_\mu+W_1)h=0$ and let $d\nu=h^2d\mu$. Let $g>0$ be smooth satisfying $(\Delta_\nu+W_2)g=0$ and let $d\theta=g^2d\nu$. Then the following are equivalent:

\noindent(a) $\Delta_\mu+W_1+W_2$ has a Green's function.

\noindent(b) $\Delta_\nu+W_2$ has a Green's function.

\noindent(c) $(M, d\theta)$ is non-parabolic.
\end{lemma}

\noindent\textbf{Proof:} Follows directly from integrating Lemma \ref{iter} over $t>0$.\qed\\

\noindent\textbf{Remark:} In Lemma \ref{itergreen}, (a) $\iff$ (b) holds even without assuming a profile $g>0$ for $\Delta_\nu+W_2$ exists. In fact, the existence of a Green's function implies the existence of a profile.

\subsection{Compactly Supported Potentials}

Let $(M, d\mu)$ be a complete, non-compact, parabolic weighted manifold satisfying (HKE) as well as (RCA) for a distinguished point $o\in M$. Let $q\geq 0$ be a smooth and compactly supported function on $M$ that is not identically zero. In this section we develop heat kernel estimates for $\Delta_\mu+q$. These are already known; see Grigor'yan \cite{grigoryan2} in particular. They will, however, be very useful for proving our main result on quadratically decaying potentials.\\

Let $h>0$ be a profile for $\Delta_\mu+q$. Let $d\nu=h^2d\mu$ and consider the weighted manifold $(M, d\nu)$. To show that $(M, d\nu)$ satisfies (HKE), it suffices to prove that $(M, d\nu)$ satisfies volume doubling and the Poincar\'e inequality. Showing these will involve the notion of remote and anchored balls. We only need a specific version of these concepts which we define below.

\begin{definition}
Given a distinguished point $o\in M$, a ball $B(x,r)$ is \textbf{remote} if $r\leq \frac{1}{2}d(o, x)$. A ball $B(x,r)$ is \textbf{anchored} if $x=o$.
\end{definition}

We also introduce the volume comparison property, abbreviated (VC).

\begin{definition}
A weighted manifold $(M, d\mu)$ satisfies the \textbf{volume comparison condition (VC)} with respect to $o\in M$ if there exists a constant $C>0$ such that for all $x\in M$ with $d(o, x)=\rho$, we have \begin{equation}
V_\mu(o, \rho)\leq CV_\mu\left(x, \frac{1}{64}\rho\right).
\end{equation}
\end{definition}

Here we cite the main theorem from Grigor'yan and Saloff-Coste \cite{grigoryansaloffcoste2}.

\begin{theorem}\label{volumecomptheorem}
(Grigor'yan, Saloff-Coste 2005) Let $(M, d\mu)$ be a complete non-compact weighted manifold satisfying (RCA) with respect to $o\in M$. Assume that $(M, d\mu)$ satisfies (VD) and (PI) for remote balls with respect to $o$. Then $(M, d\mu)$ satisfies (VD) and (PI) for all balls if and only if it satisfies (VC) with respect to $o$.
\end{theorem}

Thus to prove that $(M, d\nu)$ satisfies (VD) and (PI), it suffices to check these properties for remote balls as well as verify (VC) with respect to $o$.

\begin{lemma} \label{hvol}
Let $(M, d\mu)$ be a complete, non-compact, parabolic weighted manifold satisfying (HKE) as well as (RCA) with respect to $o\in M$. Let $q\in C_c^\infty(M)$ be nonnegative and not identically zero, let $h>0$ be a profile for $\Delta_\mu+q$, and let $d\nu=h^2d\mu$. Let $V_\nu(x,r)=\nu(B(x,r))$. We have that \begin{equation}
V_\nu(x,r)\asymp V_\mu(x,r)(H(\abs{x})+H(r))^2
\end{equation} for all $x\in M$ and all $r>0$.
\end{lemma}

\noindent\textbf{Proof:} the same proof as Lemma 4.8 in Grigor'yan and Saloff-Coste \cite{grigoryansaloffcoste1}.\qed

\begin{theorem} \label{hphi}
Let $(M, d\mu)$ be a complete, non-compact, parabolic weighted manifold satisfying (HKE) as well as (RCA) with respect to $o\in M$. Let $q\in C_c^\infty(M)$ be nonnegative and not identically zero, let $h>0$ be a profile for $\Delta_\mu+q$, and let $d\nu=h^2d\mu$. The weighted manifold $(M, d\nu)$ satisfies (VD) and (PI), and hence also (HKE).
\end{theorem}

\textbf{Proof:} $(M, h^2d\mu)$ satisfies (VD) and (PI) for all remote balls by Lemma 4.7 in Grigor'yan and Saloff-Coste \cite{grigoryansaloffcoste1}. Thus we must check that it satisfies (VC). Let $x\in M$ and let $\rho=d(x,o)$. We have that $V_\nu(o, \rho)=\int_{B(o, r)}h^2d\mu\leq h(r)^2V(o,r)\preceq h(r)^2V(x,r)\asymp V_\nu\left(x,\frac{1}{64}\rho\right)$.\qed\\



From here we could then write down a precise estimate for the heat kernel of $\Delta_\mu+q$, but we do not need it right now.\\

Note that in the theorem above, $(M, h^2d\mu)$ is non-parabolic. One way to see this is that, letting $d\nu=h^2d\mu$, $(M, d\nu)$ has a non-constant bounded $\Delta_\nu$-superharmonic function, namely $1/h$. Therefore it has a Green's function, which we will denote $G_\nu(x,y)$. We now turn to estimating $G_\nu(x,y)$.

\begin{lemma}\label{greenupper}
There exists a constant $C>0$ such that for all $x,y\in M$ distinct, we have \begin{equation}
G_\nu(x,y)\leq \frac{C}{h(y)^2}\left(\int_{d(x,y)^2}^{\max(1, \abs{y}^2)}\frac{dt}{V_\mu(y,\sqrt{t})}\right)_++\frac{C}{h(y)}.
\end{equation}
\end{lemma}

\noindent\textbf{Proof:} Lemma \ref{greenhke} tells us that (HKE) implies the following estimate on the Green's function:

\begin{equation}
G_\nu(x,y)\asymp \int_{d(x,y)^2}^\infty \frac{dt}{V_\nu(y, \sqrt{t})}\asymp \int_{d(x,y)^2}^\infty \frac{dt}{(H(\abs{y})+H(\sqrt{t}))^2V_\mu(y, \sqrt{t})}.
\end{equation}
 Note that $H'(r)=\frac{re^{-1/r}}{V(r)}$ and by the chain rule, if $f(t)=H(\abs{y})+H(\sqrt{t})$ then $f'(t)=\frac{e^{-1/\sqrt{t}}}{2V(\sqrt{t})}.$ In particular, for $t>1$ we have $f'(t)\asymp V_\mu(\sqrt{t})^{-1}$. Also, if $t\geq\abs{y}^2$, then $V_\mu(y, \sqrt{t})\asymp V_\mu(\sqrt{t})$ by volume doubling. Thus if $d(x,y)\geq \max(1, \abs{y})$, we have that $$\int_{d(x,y)^2}^\infty \frac{dt}{(H(\abs{y})+H(\sqrt{t}))^2V_\mu(y,\sqrt{t})}\asymp \int_{d(x,y)^2}^\infty \frac{f'(t)}{f(t)^2}dt= \frac{1}{H(\abs{y})+H(d(x,y))}.$$ By the same token, $$\int_{\max(1, \abs{y}^2)}^\infty\frac{dt}{(H(\abs{y})+H(\sqrt{t}))^2V_\mu(y,\sqrt{t})} \asymp \frac{1}{H(\abs{y})+H(\max(1, \abs{y}))}\asymp \frac{1}{H(\abs{y})}.$$ Thus we have \begin{equation}\label{greenbestest}G_\nu(x,y)\asymp \frac{1}{(H(\abs{y})+H(d(x,y)))^2}\left(\int_{d(x,y)^2}^{\max(1, \abs{y}^2)}\frac{dt}{V_\mu(y,\sqrt{t})}\right)_+ +\frac{1}{H(\abs{y})+H(d(x,y))},\end{equation} so in particular $$G_\nu(x,y)\preceq \frac{1}{h(y)^2}\left(\int_{d(x,y)^2}^{\max(1, \abs{y}^2)}\frac{dt}{V_\mu(y,\sqrt{t})}\right)_++\frac{1}{h(y)}.$$\qed\\

\noindent\textbf{Remark:} Under the same hypotheses as Lemma \ref{greenupper}, one has the following two-sided on $G_\nu(x,y)$: $$G_\nu(x,y)\asymp \frac{1}{(H(\abs{y})+H(\abs{x}))^2}\left(\int_{d(x,y)^2}^{\max(1, \abs{y}^2, \abs{x}^2)}\frac{dt}{V_\mu(y, \sqrt{t})}\right)_++\frac{1}{H(\abs{y})+H(\abs{x})}.$$ This is symmetric in $x$ and $y$ due to volume doubling.


\section{Gaugeability in Non-Parabolic Manifolds}\label{gaugesection}

In this section we move from considering compactly supported potentials to more general ones. In the course of applying the $h$-transform in the previous section, we turned a parabolic weighted manifold into a non-parabolic one. Working on a non-parabolic manifold will afford us various useful notions of boundedness with respect to the Green's function. Let $(M, d\mu)$ be a complete, non-compact, non-parabolic weighted manifold. Let $\Delta_\mu$ be the Laplacian of $(M, d\mu)$ and assume $\Delta_\mu\geq 0$. Let $G_\mu(x,y)$ be the Green's function of $(M, d\mu)$ and let $(X_t)_{t>0}$ denote Brownian motion on $(M, d\mu)$. \\

\noindent\textbf{Remark:} $(M, d\mu)$ being non-parabolic is equivalent to Brownian motion $(X_t)_{t>0}$ being transient on $M$.\\

The Markov transition kernel of Brownian motion on $(M, d\mu)$ is the heat kernel $p_\mu(t,x,y)$ of $\Delta_\mu$. We let $E^x$ denote expectation with respect to Brownian motion started at $x\in M$. Let $\fct{W}{M}{\Real}$ be smooth.

\begin{definition}\label{greenb}
Let $(M, d\mu)$ be a complete, non-compact, non-parabolic weighted manifold with Green's function $G_\mu(x,y)$. $W\in C^\infty(M)$ is \textbf{Green-bounded on $(M, d\mu)$} if $$\sup_{x\in M} \int_M G_\mu(x,y)\abs{W(y)}d\mu(y)<+\infty.$$
\end{definition}
When $W$ is Green-bounded we let $\norm{W}_{K^\infty}=\sup_{x\in M} \int_M G_\mu(x,y)\abs{W(y)}d\mu(y).$ We now define the Kato class at infinity.


\begin{definition}\label{kinf}
Let $(M, d\mu)$ be a complete, non-compact, non-parabolic weighted manifold. We say that $W\in C^\infty(M)$ is in the \textbf{Kato infinity class} of $(M, d\mu)$, $W\in K^\infty(M, d\mu)$, if for some (any) exhaustion $\set{\Omega_k}_{k=0}^\infty$ of smooth, relatively compact domains such that $\Omega_0\neq \emptyset$, $\overline{\Omega_k}\subseteq \Omega_{k+1}$ and $M=\bigcup_{k=0}^\infty \Omega_k$, letting $\Omega_k^*=M\setminus \Omega_k$ we have $$\lim_{k\rightarrow +\infty} \sup_{x\in M}\int_{\Omega_k^*} G_\mu(x,y)\abs{W(y)}d\mu(y)=0.$$
\end{definition}

It is routine to check that this definition is independent of the choice of the exhaustion $\set{\Omega_k}_{k=0}^\infty$, and that $W\in K^\infty(M, d\mu)$ implies $W$ is Green-bounded on $(M, d\mu)$. We have the following lemma, which is easy to check.

\begin{lemma}\label{uniformint}
Let $W\in K^\infty(M, d\mu)$. Then the family $$\set{\abs{W(y)}G_\mu(x,y)\colon x\in M}$$ is uniformly integrable on $(M, d\mu)$.
\end{lemma}

The following lemma is a setup for results in Section 5. A similar result in Euclidean space appears in Section 8 in Pinchover \cite{pinchover4}.

\begin{lemma}\label{3glemma}
Let $(M, d\mu)$ be a complete, non-compact, non-parabolic weighted manifold satisfying (HKE). Let $W\in K^\infty(M, d\mu)$ be smooth and let $\set{\Omega_k}_{k=0}^\infty$ be an exhaustion of $M$ as in Definition \ref{kinf}. For each $k$ let $\Omega_k^*=M\setminus \Omega_k$. Then $$\lim_{k\rightarrow\infty}\sup_{x,y\in \Omega_k^*}\int_{\Omega_k^*} \frac{G_\mu(x,z)G_\mu(z,y)\abs{W(z)}}{G_\mu(x,y)}d\mu(z)=0.$$

\end{lemma}

\noindent\textbf{Proof:} Let us start with a claim known as the $3G$ Principle.

\noindent\textbf{Claim:} Under the hypotheses of Lemma \ref{3glemma}, there exists $C>0$ such that for all $x,y,z\in M$ distinct, $$\frac{G_\mu(x,z)G_\mu(z,y)}{G_\mu(x,y)}\leq C(G_\mu(x,z)+G_\mu(z,y)).$$

\noindent\textbf{Proof of Claim:} Let $x,y,z\in M$ be distinct. Let $d(\cdot, \cdot)$ be the geodesic distance function. By the triangle inequality, $d(x,y)\leq d(x,z)+d(z,y)$ and therefore either $d(x,y)\leq 2d(x,z)$ or $d(x,y)\leq 2d(z,y)$. Assume that $d(x,y)\leq 2d(x,z)$. Using Lemma \ref{greenhke}, we have that $$G_\mu(x,y)\asymp \int_{d(x,y)^2}^\infty \frac{dt}{V_\mu(x, \sqrt{t})}\geq \int_{4d(x,z)^2}^\infty \frac{dt}{V_\mu(x,\sqrt{t})}\asymp G_\mu(x,z),$$ the last relation obtained by the using change of variables $s=\sqrt{2}t$ and applying volume doubling. Therefore in this case we have $G_\mu(x,z)\preceq G_\mu(x,y)$ and so $$\frac{G_\mu(x,z)G_\mu(z,y)}{G_\mu(x,y)}\preceq G_\mu(z,y).$$ If instead $d(x,y)\leq 2d(z,y)$, the proof is the same except to first use that $G_\mu(x,y)=G_\mu(y,x)$. This proves the claim. \qed\\

Now assume that $W\in K^\infty(M, d\mu)$ and let $\set{\Omega_k}_{k=0}^\infty$ be an exhaustion of $M$ as in Definition \ref{kinf}. Let $\epsilon>0$. By definition, there exists $k_0\in \N$ such that for $k\geq k_0$, $$\sup_{x\in \Omega_k^*} \int_{\Omega_k^*}G_\mu(x,y)\abs{W(y)}d\mu(y)<\frac{\epsilon}{2C}.$$ Here $C>0$ is as in the claim. Let $k\geq k_0$ and let $x,y\in \Omega_k^*$. We have $$\int_{\Omega_k^*} \frac{G_\mu(x,z)G_\mu(z,y)\abs{W(z)}}{G_\mu(x,y)}d\mu(z)\leq C\int_{\Omega_k^*} \left(G_\mu(x,z)+G_\mu(z,y)\right)\abs{W(z)}d\mu(z)< \epsilon.$$ This proves the lemma.\qed\\

\noindent\textbf{Remark:} Potentials $W$ satisfying the conclusion of Lemma \ref{3glemma} are known in the literature as \textit{small perturbations}. The ratio of three Green's functions that appears in Lemma \ref{3glemma} will be used later in Lemma \ref{pinchoverlemma}. It is also related to conditional Brownian motion; see for example Zhao \cite{zhao1} for more details.


\begin{lemma}\label{greenexpect}
Let $(M, d\mu)$ be a complete, non-compact, non-parabolic weighted manifold. Let $\fct{W}{M}{\Real}$ be smooth. Then $W$ is Green-bounded on $(M, d\mu)$ if and only if $$\sup_{x\in M}E^x\left[\int_0^\infty \abs{W(X_s)}ds\right]<+\infty.$$
\end{lemma}

\noindent\textbf{Proof:} Note that by the Fubini-Tonelli theorem and the fact that $p_\mu(t,x,y)$ is the transition kernel for $(X_t)_{t>0}$, $$\int_M G_\mu(x,y)\abs{W(y)}d\mu(y)=\int_M \left(\int_0^\infty p_\mu(t,x,y)dt\right)\abs{W(y)}d\mu(y)$$ $$=\int_0^\infty \int_M p_\mu(t,x,y)\abs{W(y)}d\mu(y)dt=\int_0^\infty E^x\left[\abs{W(X_t)}\right]dt=E^x\left[\int_0^\infty \abs{W(X_s)}ds\right].$$ The result follows.\qed

\begin{definition}
The \textbf{Feynman-Kac gauge} of $W$ is the function $$e_W(x):=E^x\left[\exp\left(-\int_0^\infty W(X_s)ds\right)\right].$$ We say that $W$ is \textbf{gaugeable} on $(M, d\mu)$ if $e_W<+\infty$.
\end{definition} Note that $e_W$ may not be finite, and may also be zero. For more on the gauge, see the book by Chung and Zhao \cite{chungzhao} as well as the articles by Chen \cite{chen1} and Zhao \cite{zhao1}. We have the following well-known result, known as Khasmin'skii's Lemma.

\begin{lemma}\label{khas} (Khasmin'skii's Lemma)
Let $W$ be such that $\norm{W}_{K^\infty}=\alpha<1$. Then $\sup_{x\in M}e_W(x)<(1-\alpha)^{-1}$.
\end{lemma}

The following is a essentially a consequence of Theorem 2 in Zhao \cite{zhao1} in Euclidean space; we reproduce the proof on a weighted manifold.

\begin{theorem}\label{gaugepro}
Let $(M, d\mu)$ be a complete, non-compact, non-parabolic weighted manifold. Let $W\in K^\infty(M, d\mu)$. Then $W$ is gaugeable on $(M, d\mu)$ if and only if there exists a smooth function $h>0$ satisfying $(\Delta_\mu+W)h=0$ and $h\asymp 1$.
\end{theorem}

\noindent\textbf{Proof:} ($\Rightarrow$) Set $h(x)=e_W(x)$. By hypothesis, $h$ is bounded above. Let $C>0$ such that $\sup_{x\in M} \int_M G_\mu(x,y)\abs{W(y)}d\mu(y)\leq C$. To see that $h$ is bounded below, note that by Jensen's inequality, $$h(x)=E^x\left[\exp\left(-\int_0^\infty W(X_s)ds\right)\right] \geq \exp\left(-E^x\left[\int_0^\infty W(X_s)ds\right]\right)\geq e^{-C},$$ the last by Lemma \ref{greenexpect}. Thus $h\asymp 1$.\\

Now we show that $(\Delta_\mu+W)h=0$. By the Fubini-Tonelli theorem as well as the Markov property and integrating by parts, we have $$\int_M G_\mu(x,y)W(y)h(y)d\mu(y)=E^x\left[\int_0^\infty W(X_t) E^{X_t}\left[\exp\left(-\int_0^\infty W(X_s)ds\right)\right]\right]$$ $$=E^x\left[\int_0^\infty W(X_t)\exp\left(-\int_t^\infty W(X_s)ds\right)\right] = E^x\left[1-\exp\left(-\int_0^\infty W(X_s)ds\right)\right]$$ $$=1-h(x).$$ Applying the distributional version of $\Delta_\mu$ to both sides of \begin{equation}\label{integralidentity}\int_M G_\mu(x,y)W(y)h(y)d\mu(y)=1-h(y)\end{equation} yields that $(\Delta_\mu+W)h=0$ as a distribution. $h$ is continuous using the result of Lemma \ref{uniformint} that $\set{G_\mu(x,y)W(y)h(y)\colon x\in M}$ is uniformly integrable. By elliptic regularity, $h$ is smooth. \\

($\Leftarrow$) Let $h>0$ be smooth such that $(\Delta_\mu+W)h=0$ and $h\asymp 1$. Let $\set{\Omega_k}_{k=0}^\infty$ be an exhaustion as in Definition \ref{kinf}. Fix $x\in M$ and let $k_0\in \N$ be such that $x\in \Omega_{k_0}$. By the Feynman-Kac formula, for all $k\geq k_0$ we have $$h(x)=E^x\left[\exp\left(-\int_0^{\tau_k} W(X_s)ds\right)h(X_{\tau_k})\right]\asymp E^x\left[\exp\left(-\int_0^{\tau_k} W(X_s)ds\right)\right].$$ Here $\tau_k=\inf\set{t\colon X_t\notin \Omega_k}$ is the first exit time of Brownian motion from $\Omega_k$. By Fatou's Lemma, we have $$E^x\left[\exp\left(-\int_0^\infty W(X_s)ds\right)\right]\leq \liminf_{k\rightarrow \infty} E^x\left[\exp\left(-\int_0^{\tau_k} W(X_s)ds\right)\right]\preceq h(x)\preceq 1.$$ Thus $W$ is gaugeable.\qed\\

\noindent\textbf{Remark:} See Murata \cite{murata5} for a version of Equation \ref{integralidentity} for more general elliptic operators.\\

We finish this section with the following lemma.

\begin{lemma}\label{greensmallgauge}
Let $W\in K^\infty(M, d\mu)$ such that $$\sup_{x\in M}\int_M G_\mu(x,y)W_-(y)d\mu(y)<1.$$ Then $\Delta_\mu+W$ is gaugeable on $(M, d\mu)$.
\end{lemma}

\noindent\textbf{Proof:} By the hypothesis and Khasmin'skii's Lemma, $\Delta_\mu-W_-$ is gaugeable on $(M, d\mu)$. By monotonicity of the Feynman-Kac gauge, $\Delta_\mu+W$ is gaugeable on $(M, d\mu)$.\qed

\section{Subcriticality and Criticality Theory}\label{criticalsection}

\begin{definition}\label{critdef}
Let $(M, d\mu)$ be a complete, non-compact weighted manifold and let $W\in C^\infty(M)$. We say that:

\noindent(i) $\Delta_\mu+W$ is \textbf{subcritical} on $(M, d\mu)$ if $\Delta_\mu+W\geq 0$ and $\Delta_\mu+W$ admits a positive Green's function.
 
\noindent (ii)  $\Delta_\mu+W$ is \textbf{critical} on $(M, d\mu)$ if $\Delta_\mu+W\geq 0$ and $\Delta_\mu+W$ is not subcritical. 

\noindent (iii) $\Delta_\mu+W$ is \textbf{supercritical} on $(M, d\mu)$ if $\Delta_\mu+W\not\geq 0$.
\end{definition}

By the Allegretto-Piepenbrink Theorem, $\Delta_\mu+W\geq 0$ if and only if $\Delta_\mu+W$ has a profile $h>0$ on $M$ satisfying $(\Delta_\mu+W)h=0$. Letting $d\nu=h^2d\mu$, we see that the Green's function $G_\mu^W(x,y)$ of $\Delta_\mu+W$, if it exists, satisfies $G^W_\mu(x,y)=h(x)h(y)G_\nu(x,y)$, where $G_\nu(x,y)$ is the Green's function of $\Delta_\nu$ on $(M, d\nu)$.

\begin{lemma} \label{hpreservescrit}
Let $(M, d\mu)$ be a complete, non-compact weighted manifold. Let $q\in C^\infty(M)$ be such that there exists a profile $h>0$ satisfying $(\Delta_\mu+q)h=0$ and let $d\nu=h^2d\mu$. Let $W\in C^\infty(M)$. Then $\Delta_\mu+W$ is subcritical on $(M, d\mu)$ if and only if $\Delta_\nu+W-q$ is subcritical on $(M, d\nu)$. The previous sentence also holds with ``critical" or ``supercritical" replacing ``subcritical."
\end{lemma}

\noindent\textbf{Proof:} By Lemma \ref{itergreen} and the subsequent remark, it suffices to prove that $\Delta_\mu+W\geq 0$ if and only if $\Delta_\nu+W-q\geq 0$. This is a routine exercise.\qed





\begin{lemma}\label{compactisgbounded}
Let $(M, d\mu)$ be a complete, non-compact, non-parabolic weighted manifold. Let $q\in C_c^\infty(M)$. Then $q$ is Green-bounded on $(M, d\mu)$.
\end{lemma}

\noindent\textbf{Proof:} Follows from Lemma 10.3 (a) in Grigor'yan \cite{grigoryan2}.\qed

\begin{lemma}\label{perturb}
Let $(M, d\mu)$ be a complete, non-compact, non-parabolic weighted manifold. Let $q\in C_c^\infty(M)$. Then there exists $c>0$ such that $\Delta_\mu+cq\geq 0$.
\end{lemma}

\noindent\textbf{Proof:} Since $q$ is Green-bounded on $(M, d\mu)$ by the previous lemma, there exists $c>0$ such that $$\sup_{x\in M}\int_M G(x,y)\abs{cq(y)}d\mu(y)<1.$$ By Theorem \ref{gaugepro}, $\Delta_\mu+cq$ has a profile $h>0$ with $(\Delta_\mu+cq)h=0$ and $h\asymp 1$. By the Allegretto-Piepenbrink Theorem, $\Delta_\mu+cq\geq 0$.\qed\\

The next lemma, which goes back to Murata and Pinchover, connects subcriticality to compactly supported perturbations. See Devyver \cite{devyver1}, Theorem 3.7 for a similar result with different proof.

\begin{lemma}\label{subcritperturb}
Let $(M, d\mu)$ be a complete, non-compact weighted manifold. Let $W\in C^\infty(M)$ be bounded. Then $\Delta_\mu+W$ is subcritical if and only if there exists $q\in C_c^\infty(M)$ non-negative and not identically zero such that $\Delta_\mu+W-q\geq 0$.
\end{lemma}

\noindent\textbf{Proof:} ($\Rightarrow$) Assume $\Delta_\mu+W$ is subcritical. Let $h>0$ be a profile for $\Delta_\mu+W$. Let $d\nu=h^2d\mu$ and consider the weighted manifold $(M, d\nu)$. By Lemma \ref{itergreen}, $\Delta_\nu$ is subcritical, i.e. $(M, d\nu)$ is non-parabolic. Let $\phi\in C_c^\infty(M)$ be non-positive and not identically zero. By Lemma \ref{perturb}, there exists $c>0$ such that $\Delta_\nu+c\phi\geq 0$. Set $q=-c\phi$. Then since $\Delta_\nu-q\geq 0$ we have $\Delta_\mu+W-q\geq 0$ as desired.\\

($\Leftarrow$) Let $q\in C_c^\infty(M)$ be non-negative and not identically zero such that $\Delta_\mu+W-q\geq 0$. By the Allegretto-Piepenbrink Theorem, there exists a profile $h>0$ for $\Delta_\mu+W-q$. Let $d\nu=h^2d\mu$. Since $\Delta_\nu\geq 0$, we clearly have $\Delta_\nu+q\geq 0$. Thus by the Allegretto-Piepenbrink Theorem, there exists a profile $g>0$ for $\Delta_\nu+q$. Let $d\theta=g^2d\nu$.\\

Now we have $$\Delta_\theta\left(\frac{1}{g}\right)=\frac{1}{g}\Delta_\nu(1)+\frac{q}{g}=\frac{q}{g}\geq 0.$$ Thus $1/g$ is $\Delta_\theta$-superharmonic. Since $1/g$ is positive and non-constant, by Liouville's Theorem \ref{liouville}, $(M, d\theta)$ is non-parabolic. Therefore $\Delta_\theta$ is subcritical. Unraveling we find via Lemma \ref{itergreen} that $\Delta_\nu+q$ is subcritical and so is $\Delta_\mu+(W-q)+q=\Delta_\mu+W$.\qed\\

We now have the following result from Devyver \cite{devyver1}.

\begin{theorem}\label{devyvertheorem} (Devyver \cite{devyver1}, Theorem 3.2) Let $(M, d\mu)$ be a complete, non-compact, non-parabolic weighted manifold. Let $W\in C^\infty(M)$ be Green-bounded such that $W_-\in K^\infty(M, d\mu)$ and such that $\Delta_\mu+W$ is subcritical. Then there exists a profile $h>0$ satisfying $(\Delta_\mu+W)h=0$ as well as $h\asymp 1$.
\end{theorem}

So subcritical operators have profiles that are approximately constant. We have the following result that ties together several concepts discussed so far. This is generally known; see Devyver \cite{devyver1} Theorem 3.2 and Theorem 3.7 as well as Zhao \cite{zhao1} Theorem 2, but we find it helpful to collect them in one place. 

\begin{theorem}\label{zhaomain}
Let $(M, d\mu)$ be a complete, non-compact, non-parabolic weighted manifold satisfying (HKE). Let $W\in K^\infty(M, d\mu)$. The following are equivalent:

\noindent(a) $\Delta_\mu+W$ is subcritical.

\noindent(b) $W$ is gaugeable on $(M, d\mu)$.

\noindent(c) There exists a profile $h>0$ for $\Delta_\mu+W$ satisfying $h\asymp 1$.

\noindent(d) There exists $q\in C_c^\infty(M)$ non-negative and not identically zero such that $\Delta_\mu+W-q\geq 0$.
\end{theorem}

\noindent\textbf{Proof:} (b) $\iff$ (c), (a) $\Rightarrow$ (c), and (a) $\iff$ (d) were already shown in Theorem \ref{gaugepro}, Theorem \ref{devyvertheorem}, and Lemma \ref{subcritperturb} respectively. Thus it suffices to show that (c) $\Rightarrow$ (a).\\


\noindent(c) $\Rightarrow$ (a): Let $h>0$ be a profile for $\Delta_\mu+W$ satisfying $h\asymp 1$ and let $d\nu=h^2d\mu$. Notably, since the property (HKE) is equivalent to the properties (VD) and (PI), it is routine to check that (VD) and (PI) for $(M, d\mu)$ imply the same for $(M, d\nu)$ and so $(M, d\nu)$ satisfies (HKE). Furthermore, $(M, d\nu)$ is non-parabolic as well and so has a Green's function $G_\nu(x,y)$. We then see that a Green's function for $\Delta_\mu+W$ exists with $G_\mu^W(x,y)=h(x)h(y)G_\nu(x,y)$.\qed\\





The following well-known proposition produces many examples of subcritical Schr\"odinger operators.

\begin{proposition}\label{positivesubcritical}
Let $(M, d\mu)$ be a complete, non-compact, weighted manifold. Let $W\in C^\infty(M)$ be bounded, nonnegative and not identically zero. Then $\Delta_\mu+W$ is subcritical on $(M, d\mu)$.
\end{proposition}

\noindent\textbf{Proof:} By hypothesis there exists $q\in C_c^\infty(M)$ with $0\leq q\leq W$. Thus $W-q\geq 0$ and so $\Delta_\mu+W-q\geq 0$. By Lemma \ref{subcritperturb}, $\Delta_\mu+W$ is subcritical.\qed\\

 All this is in contrast to the critical case, which we now cite from Pinchover \cite{pinchover1}.

\begin{lemma} \label{pinchoverlemma} (Pinchover \cite{pinchover1}, Lemma 2.5) Let $(M, d\mu)$ be a complete, non-compact, non-parabolic weighted manifold. Let $\set{\Omega_k}_{k=0}^\infty$ be an exhaustion of $M$ as in Definition \ref{kinf}, and let $\Omega_k^*=M\setminus \Omega_k$. Let $W\in C^\infty(M)$ and assume that $\Delta_\mu+W$ is critical and that $$\lim_{k \rightarrow\infty} \sup_{x,y\in \Omega_k^*}\int_{\Omega_k^*} \frac{G_\mu(x,z)G_\mu(z,y)\abs{W(z)}}{G_\mu(x,y)}d\mu(z)=0.$$ Let $o\in M$ be a distinguished point and let $h>0$ be a profile for $\Delta_\mu+W$. Outside a compact set containing $o$ we have $$h(x)\asymp G_\mu(x,o).$$

\end{lemma}

With further estimates on the Green's function, we can estimate the profile of critical Schr\"odinger operator.

\begin{lemma}\label{criticalprofile}
Let $(M, d\mu)$ be a complete, non-compact, parabolic weighted manifold satisfying (HKE) as well as (RCA) with respect to a distinguished point $o\in M$. Let $q\in C_c^\infty(M)$ be nonnegative and not identically zero, let $h>0$ be a profile for $\Delta_\mu+q$, and let $d\nu=h^2d\mu$. Let $W\in K^\infty(M, d\nu)$ be such that $\Delta_\nu+W$ is critical, and let $g>0$ be a profile for $\Delta_\nu+W$. Then $$g\asymp \frac{1}{h}.$$
\end{lemma}

\noindent\textbf{Proof:} Because $(M, d\nu)$ is non-parabolic and satisfies (HKE), Lemma \ref{3glemma} implies that the hypotheses of Lemma \ref{pinchoverlemma} are satisfied. Therefore by Lemma \ref{pinchoverlemma} $$g(x)\asymp G_\nu(x,o)$$ outside a compact set. However, by Expression \ref{greenbestest} in the proof of Lemma \ref{greenupper} we have, for $\abs{x}\geq 1$, $$G_\nu(x,o)=G_\nu(o,x)\asymp \frac{1}{(H(\abs{x})+H(\abs{x}))^2}\left(\int_{\abs{x}^2}^{\max(1, \abs{x}^2)}\frac{dt}{V_\mu(x,\sqrt{t})}\right)_+ +\frac{1}{H(\abs{x})+H(\abs{x})}$$ $$\asymp \frac{1}{H(\abs{x})}\asymp \frac{1}{h(x)}.$$ Therefore $$g\asymp \frac{1}{h}$$ outside of a compact set, which can be extended to all of $M$ by noting that each function is approximately a constant on any fixed compact set.\qed\\

\noindent\textbf{Remark:} The results in this section when $M=\Real^2$ should be compared to those of Murata in \cite{murata2} and \cite{murata4}, where similar results were obtained for $W$ satisfying a form of polynomial decay.

\section{Potentials with Faster Than Quadratic Decay: Main Results}\label{mainresultssection}

We now turn to determining which potentials are Green-bounded and Kato infinity class.

\begin{lemma}\label{gbound}
Let $(M, d\mu)$ be a complete, non-compact, parabolic weighted manifold satisfying (HKE) as well as (RCA) with respect to $o\in M$. Let $$H(\abs{y})=1+\left(\int_1^{\abs{y}^2} \frac{dt}{V_\mu(\sqrt{t})}\right)_+$$ and let $$\widehat{H}(y)=1+\left(\int_1^{\abs{y}^2}\frac{dt}{V_\mu(y, \sqrt{t})}\right)_+.$$ 
Let $q\in C_c^\infty(M)$ be nonnegative and not identically zero, let $h>0$ be a profile for $\Delta_\mu+q$, and let $d\nu=h^2d\mu$. Let $W\in C^\infty(M)$ be bounded such that $$\int_M \abs{W(y)}(H(\abs{y})+\widehat{H}(y))d\mu(y)<+\infty.$$ Then $W$ is Green-bounded on $(M, d\nu)$.
\end{lemma}

\noindent\textbf{Proof:} Let $W$ be as in the hypothesis, let $C=\sup_M\abs{W}$ and let $I=\int_M \abs{W(y)}h(y)d\mu(y)$. We have that $$\int_M G_\nu(x,y)\abs{W(y)}d\nu(y)=\int_M G_\nu(x,y)\abs{W(y)}h(y)^2d\mu(y)$$ $$\preceq \int_M\abs{W(y)}\left(\int_{d(x,y)^2}^{\max(1, \abs{y}^2)}\frac{dt}{V_\mu(y,\sqrt{t})}\right)_+d\mu(y)+\int_M\abs{W(y)}h(y)d\mu(y).$$

Let us estimate $\int_M\abs{W(y)}\left(\int_{d(x,y)^2}^{\max(1, \abs{y}^2)}\frac{dt}{V_\mu(y,\sqrt{t})}\right)_+d\mu(y)$ by breaking the domain of integration into pieces, namely $M=B(x, 1)\cup (M\setminus B(x, 1))$. \\

We have that $$\int_{B(x, 1)} \abs{W(y)}\left(\int_{d(x,y)^2}^{\max(1,\abs{y}^2)} \frac{dt}{V_\mu(y,\sqrt{t})}\right)_+d\mu(y)$$ $$\leq \int_{B(x, 1)} \abs{W(y)} \left(\int_{d(x,y)^2}^{1} \frac{dt}{V_\mu(y,\sqrt{t})}\right)_+d\mu(y)+\int_{B(x, 1)\setminus B(o, 1)} \abs{W(y)} \left(\int_{1}^{\abs{y}^2} \frac{dt}{V_\mu(y,\sqrt{t})}\right)d\mu(y)$$ $$\leq \left(\sup_{B(x, 1)}\abs{W}\right)\int_{B(x, 1)}\left(\int_{d(x,y)^2}^{1}\frac{dt}{V_\mu(y,\sqrt{t})}\right)_+d\mu(y)+\int_{B(x, 1)}\abs{W(y)}\widehat{H}(y)d\mu(y).$$\\

 Note that if $y\in B(x, \sqrt{t})$, then $B(x, \sqrt{t})\subseteq B(y, 2\sqrt{t})$. Thus by volume doubling, $V_\mu(x, \sqrt{t})\leq C_{VD}V_\mu(y, \sqrt{t})$. By the Fubini-Tonelli theorem, we have $$ \int_{B(x, 1)}\left(\int_{d(x,y)^2}^{1}\frac{dt}{V_\mu(y,\sqrt{t})}\right)_+d\mu(y)=\int_0^{1} \int_{B(x, \sqrt{t})}\frac{1}{V_\mu(y, \sqrt{t})}d\mu(y)dt$$ $$\leq \int_0^1 \int_{B(x, \sqrt{t})}\frac{C_{VD}}{V_\mu(x, \sqrt{t})}d\mu(y)dt\leq C_{VD}.$$\\

The next piece is integrating over $M\setminus B(x, 1)$. We have that $$\int_{M\setminus B(x, 1)} \abs{W(y)}\left(\int_{d(x,y)^2}^{\max(1, \abs{y}^2)}\frac{dt}{V_\mu(y,\sqrt{t})}\right)_+d\mu(y)$$ $$\leq\int_{M\setminus B(x, 1)} \abs{W(y)}\left(\int_{1}^{\abs{y}^2}\frac{dt}{V_\mu(y,\sqrt{t})}\right)_+d\mu(y)\leq \int_{M\setminus B(x, 1)} \abs{W(y)}\widehat{H}(y)d\mu(y).$$ Putting it all together we have that $$\int_M \abs{W(y)}\left(\int_{d(x,y)^2}^{\max(1, \abs{y}^2)}\frac{dt}{V_\mu(\sqrt{t})}\right)_+d\mu(y)$$  $$=\int_{B(x, 1)}\abs{W(y)}\left(\int_{d(x,y)^2}^{\max(1, \abs{y}^2)}\frac{dt}{V_\mu(\sqrt{t})}\right)_+d\mu(y)+\int_{M\setminus B(x, 1)}\abs{W(y)}\left(\int_{d(x,y)^2}^{\max(1, \abs{y}^2)}\frac{dt}{V_\mu(\sqrt{t})}\right)_+d\mu(y)$$ $$\preceq C_{VD}\sup_{B(x, 1)} \abs{W}+\int_M \abs{W(y)}\widehat{H}(y)d\mu(y).$$ Thus we have that $$\int_M G_\nu(x,y)\abs{W(y)}d\nu(y)$$ $$\preceq \int_M\abs{W(y)}\left(\int_{d(x,y)^2}^{\max(1, \abs{y}^2)}\frac{dt}{V_\mu(\sqrt{t})}\right)_+d\mu(y)+\int_M\abs{W(y)}h(y)d\mu(y)$$ $$\leq C_{VD}\sup_{B(x, 1)} \abs{W} + \int_M \abs{W(y)}(H(\abs{y})+\widehat{H}(y))d\mu(y)<+\infty.$$ Since this bound is uniform in $x\in M$, we are done.\qed

\begin{lemma}\label{winkinf}
Let $(M, d\mu)$ be a complete, non-compact, parabolic weighted manifold satisfying (HKE) as well as (RCA) with respect to $o\in M$. Let $q\in C_c^\infty(M)$ be nonnegative and not identically zero, let $h>0$ be a profile for $\Delta_\mu+q$, and let $d\nu=h^2d\mu$. Let $W\in C_0^\infty(M)$ such that $$\int_M \abs{W(y)}(H(\abs{y})+\widehat{H}(y))d\mu(y)<+\infty.$$ Then $W\in K^\infty(M, d\nu)$.
\end{lemma}

\noindent\textbf{Proof:} The same proof as that of Lemma \ref{gbound} above, except using that if $\set{\Omega_k}_{k=0}^\infty$ is an exhaustion of $M$ as in Definition \ref{kinf}, then $$\lim_{k\rightarrow \infty}\sup_{\Omega_k^*}\abs{W}=\lim_{k\rightarrow\infty} \int_{\Omega_k^*} \abs{W(y)}(H(\abs{y})+\widehat{H}(y))d\mu(y)=0.$$\qed

\begin{proposition}\label{quad}
Let $(M, d\mu)$ be a complete, non-compact, parabolic weighted manifold satisfying (HKE) as well as (RCA) with respect to $o\in M$.  Assume further there exists $C>0$ with $$\widehat{H}(y)\leq C H(\abs{y})$$ for all $y\in M$, and that there exist $\delta>0$ and $C_\delta>0$ such that for all $1\leq s<r$ we have \begin{equation}
\frac{V_\mu(r)}{V_\mu(s)}\leq C_\delta\left(\frac{r}{s}\right)^{2+\delta}.
\end{equation} Let $W\in C^\infty(M)$ be such that for constants $c, \epsilon>0$ and all $x\in M$ we have \begin{equation}\abs{W(x)}\preceq \brac{x}^{-(2+\delta+\epsilon)}.\end{equation} Let $q\geq 0$ be smooth, compactly supported and not identically zero, let $h>0$ be a profile for $\Delta_\mu+q$, and let $d\nu=h^2d\mu$. Then $W \in K^\infty(M, d\nu)$.
\end{proposition}

\noindent\textbf{Proof:} By Lemma \ref{winkinf}, it suffices to prove that $\int_M \abs{W(y)}(H(\abs{y})+\widehat{H}(y))d\mu(y)<+\infty.$ We now have $$\int_M \abs{W(y)}(H(\abs{y})+\widehat{H}(y))d\mu(y)\preceq \int_M \abs{W(y)}H(\abs{y})d\mu(y)$$ $$\preceq \int_{B(o, 1)} \abs{W(y)}H(\abs{y})d\mu(y)+\sum_{n=0}^\infty\int_{B(o, 2^{n+1})\setminus B(o, 2^n)} (2^n)^{-(2+\delta+\epsilon)}\left(\int_1^{2^n} \frac{s}{V_\mu(s)}ds\right) d\mu$$ $$\preceq 1+\sum_{n=0}^\infty(2^n)^{-(2+\delta+\epsilon)}\left(\int_1^{2^n} \frac{s}{V_\mu(s)}ds\right)V_\mu(2^n)=1+\sum_{n=0}^\infty\left(\int_1^{2^n} \frac{V_\mu(2^n)}{V_\mu(s)}sds\right)(2^n)^{-(2+\delta+\epsilon)}$$ $$\leq 1+\sum_{n=0}^\infty\left(\int_1^{2^n} C_{\delta}\left(\frac{2^n}{s}\right)^{2+\delta} sds\right)(2^n)^{-(2+\delta+\epsilon)}=1+C_{\delta}\sum_{n=0}^\infty \left(\int_1^{2^n}s^{-(1+\delta)}ds\right)(2^n)^{-\epsilon} $$ $$\leq 1+C_{\delta}\left(\frac{1}{\delta}\right)\sum_{n=0}^\infty (2^{-\epsilon})^{n}=1+C_{\delta}\left(\frac{1}{\delta}\right)\frac{1}{1-2^{-\epsilon}}<+\infty.$$ \qed\\

\noindent\textbf{Remark:} The condition $\widehat{H}(y)\leq CH(\abs{y})$ follows if we know that $V_\mu(r)\preceq V_\mu(y, r)$ over all $y\in M$, $0<r<\abs{y}$. We will see later in \ref{modelmanifolds} an example that satisfies the first of these properties but not the second.\\

This brings us to our main theorem. Let us recall once again that $H(r)= 1+\left(\int_1^{r^2} \frac{dt}{V_\mu(\sqrt{t})}\right)_+$ and $\widehat{H}(y)=1+\left(\int_1^{\abs{y}^2}\frac{dt}{V_\mu(y, \sqrt{t})}\right)_+$.






\begin{theorem} \label{maithm} Let $(M, d\mu)$ be a complete, non-compact, parabolic weighted manifold satisfying (HKE) as well as (RCA) with respect to $o\in M$. Let $W\in C_0^\infty(M)$ be such that $$\int_{M}\abs{W(y)}(H(\abs{y})+\widehat{H}(y))d\mu(y)<+\infty.$$ Let $p_\mu^W(t,x,y)$ denote the heat kernel of $\Delta_\mu+W$. There exist constants $c_1, c_2, c_3, c_4>0$ such that:

\noindent\textbf{(i)} If $\Delta_\mu+W$ is subcritical on $(M, d\mu)$, then for all $t>0$ and all $x,y\in M$ we have $$\frac{c_1H(\abs{x})H(\abs{y})}{(H(\abs{y})+H(\sqrt{t}))^2V_\mu(y,\sqrt{t})}e^{-c_2\frac{d(x,y)^2}{t}}\leq p^{W}_\mu(t,x,y)\leq\frac{c_3H(\abs{x})H(\abs{y})}{(H(\abs{y})+H(\sqrt{t}))^2V_\mu(y,\sqrt{t})}e^{-c_4\frac{d(x,y)^2}{t}}.$$

\noindent\textbf{(ii)} If $\Delta_\mu+W$ is critical on $(M, d\mu)$, then for all $t>0$ and $x,y\in M$ we have $$\frac{c_1}{V_\mu(y,\sqrt{t})}e^{-c_2\frac{d(x,y)^2}{t}}\leq p^{W}_\mu(t,x,y)\leq\frac{c_3}{V_\mu(y,\sqrt{t})}e^{-c_4\frac{d(x,y)^2}{t}}.$$
\end{theorem}

\noindent\textbf{Proof:} Let $q\in C_c^\infty(M)$ be nonnegative and not identically zero. Let $h>0$ be such that $(\Delta_\mu+q)h=0$ and let $d\nu=h^2d\mu$. By Lemma \ref{winkinf} we have that $W-q\in K^\infty(M, d\nu)$.\\

\noindent\textbf{\textit{(i)}} Since $\Delta_\mu+W$ is subcritical on $(M, d\mu)$, $\Delta_\nu+W-q$ is subcritical on $(M, d\nu)$. Hence by Theorem \ref{devyvertheorem}, there exists a profile $g>0$ for $\Delta_\nu+W-q$ with $g\asymp 1$. Thus if we let $d\theta = g^2d\nu$, the weighted manifold $(M, d\theta)$ satisfies (VD) and (PI) since $(M, d\nu)$ does. Therefore $(M, d\theta)$ satisfies (HKE). The result then follows from Lemma \ref{iter} as well as using Lemma \ref{hvol} to get $V_\theta(x,\sqrt{t})\asymp (H(\abs{y})+H(\sqrt{t}))^2V_\mu(x,\sqrt{t})$.\qed\\

\noindent\textbf{\textit{(ii)}} Next, instead assume that $\Delta_\mu+W$ is critical on $(M, d\mu)$. By Lemma \ref{hpreservescrit}, $\Delta_\nu+W-q$ is critical on $(M, d\nu)$. Let $g>0$ be a profile for $\Delta_\nu+W-q$. By Lemma \ref{criticalprofile} we have that $g\asymp 1/h$. Therefore $gh\asymp 1$ and so $(M, d\theta)$ as in the proof of \textit{(i)} satisfies (HKE) because $(M, d\mu)$ does. Furthermore, we have $V_\theta(x, \sqrt{t})\asymp V_\mu(x,\sqrt{t})$. The result follows.\qed\\

\noindent\textbf{Remark:} This result should be compared with Theorem 10.10 of Grigor'yan \cite{grigoryan2} as well as our combined statement of Murata and Grigor'yan in Theorem \ref{muratagrigoryan}. In fact, using our estimates on the profiles, Theorem \ref{maithm} can be proved using Grigor'yan's result, which is stated when $W=\frac{\Delta_\mu \phi}{\phi}$, with assumptions on the behavior of $\phi>0$. What we have supplied is an estimate of the profile $h>0$ for $\Delta_\mu+W$ when $W$ satisfies the hypothesis of Theorem \ref{maithm}; in part (i), we have $h(x)\asymp H(\abs{x})$, and in part (ii), we have $h(x)\asymp 1$. In this sense we can view our result as building on that of Murata in \cite{murata2} stated in Theorem \ref{muratatheorem}.  

\section{Applications and Examples}\label{examplesection}

By Proposition \ref{positivesubcritical}, positive smooth potentials $W$ yield subcritical operators $\Delta_\mu+W$. Let us now address a wide class of $W$ to which our Main Theorem \ref{maithm} applies: $W$ that decay to zero at infinity faster than quadratically.

\begin{theorem}\label{quaddecayapplication}
Let $(M, d\mu)$ be a complete, non-compact, parabolic weighted manifold satisfying (HKE) as well as (RCA) with respect to $o\in M$. Assume that there exists $C>0$ such that $$\widehat{H}(y)\leq CH(\abs{y})$$ for all $y\in M$ and that there exist $\delta>0$ and $C_\delta>0$ such that for all $1\leq s<r$ we have $$\frac{V_\mu(r)}{V_\mu(s)}\leq C_\delta\left(\frac{r}{s}\right)^{2+\delta}.$$ Let $W\in C^\infty(M)$ be not identically zero for which there exist $c,\epsilon>0$ such that for all $x\in M$, $$0\leq W(x)\leq c\brac{x}^{-(2+\delta+\epsilon)}.$$ Letting $p_\mu^W(t,x,y)$ denote the heat kernel of $\Delta_\mu+W$, there exist constants $c_1, c_2, c_3, c_4$ such that for all $t>0$ and $x,y\in M$ we have $$\frac{c_1H(\abs{x})H(\abs{y})}{(H(\abs{y})+H(\sqrt{t}))^2V_\mu(y,\sqrt{t})}e^{-c_2\frac{d(x,y)^2}{t}}\leq p^{W}_\mu(t,x,y)\leq\frac{c_3H(\abs{x})H(\abs{y})}{(H(\abs{y})+H(\sqrt{t}))^2V_\mu(y,\sqrt{t})}e^{-c_4\frac{d(x,y)^2}{t}}.$$
\end{theorem}

\noindent\textbf{Proof:} Since $W\geq 0$ is not identically zero, $\Delta_\mu+W$ is subcritical by Proposition \ref{positivesubcritical}. By Proposition \ref{quad}, $W$ satisfies the hypotheses of the Main Theorem \ref{maithm}(i), and the conclusion follows.\qed\\

Consider the case when $(M, d\mu)=(\Real^2, dx)$, where $\Real^2$ has the Euclidean metric and $dx$ denotes Lebesgue measure. Potentials that decay faster than quadratically as in Theorem \ref{quaddecayapplication} are already covered by the results of Murata and Grigor'yan mentioned in the introduction as Theorem \ref{muratagrigoryan}. We note, however, that the Main Theorem \ref{maithm} still offers the following addition.


\begin{corollary}\label{plane}
Let $W\in C_0^\infty(\Real^2)$ be such that $\log\brac{x}\abs{W(x)}\in L^1(\Real^2, dx)$. Let $p^W(t,x,y)$ denote the heat kernel of $\Delta+W$. There exist constants $c_1, c_2, c_3, c_4>0$ such that:

\noindent(i) If $\Delta+W$ is subcritical, then for all $t>0$, $x,y\in \Real^2$, $$\frac{c_1\log\brac{x}\log \brac{y}}{t\log(\brac{x}+\sqrt{t})\log(\brac{y}+\sqrt{t})}e^{-c_2\frac{\abs{x-y}^2}{t}}\leq p^W(t,x,y)\leq \frac{c_3 \log\brac{x}\log \brac{y}}{t\log(\brac{x}+\sqrt{t})\log(\brac{y}+\sqrt{t})}e^{-c_4\frac{\abs{x-y}^2}{t}}.$$

\noindent(ii) If $\Delta+W$ is critical, then for all $t>0$, $x,y\in \Real^2$, $$\frac{c_1}{t}e^{-c_2\frac{\abs{x-y}^2}{t}}\leq p^W(t,x,y)\leq \frac{c_3}{t}e^{-c_4\frac{\abs{x-y}^2}{t}}.$$
\end{corollary}

\noindent\textbf{Proof:} It is routine to check that $(\Real^2, dx)$ is complete, non-compact, parabolic, satisfies (HKE) as well as (RCA) with respect to the origin. If $W\in C_0^\infty(\Real^2)$ is such that $\log\brac{x}\abs{W(x)}\in L^1(\Real^2, dx)$, then by noting that $H(\abs{y})\asymp \widehat{H}(y)\asymp \log\brac{y}$ we can apply Theorem \ref{maithm} as well as symmetrization.\qed\\

\noindent\textbf{Remark:} Subcritical operators are plentiful, as any non-negative, not identically zero potential satisfying the hypothesis will yield a subcritical operator. On the other hand, critical operators other than $\Delta_\mu$ are difficult to describe, but we show that many exist in Proposition \ref{criticalexample}.

\subsection{Existence of Critical Operators}

The suspicious reader of Theorem \ref{maithm}(ii) may wonder whether, on a parabolic manifold, there exist critical Schr\"odinger operators $\Delta_\mu+W$ other than $W=0$. It is difficult to write down examples explicitly, but we can show that many exist. See Pinchover \cite{Pinchover1998} Theorem 4.1 for a related result.

\begin{proposition}\label{criticalexample}
Let $(M, d\mu)$ be a complete, non-compact, parabolic weighted manifold satisfying (HKE) as well as (RCA) with respect to $o\in M$. Let $q\in C_c^\infty(M)$ be nonnegative and not identically zero, let $h>0$ be a profile for $\Delta_\mu+q$, and let $d\nu=h^2d\mu$. Let $W_1, W_2\in C^\infty(M)$ be bounded such that:

\noindent (i) $\Delta_\mu+W_1$ is subcritical on $(M, d\mu)$,

\noindent (ii) $W_1, W_2\in K^\infty(M, d\nu)$, and

\noindent (iii) $W_2(x)> q(x)$ for some point $x\in M$.

\noindent Then there exists $c>0$ such that $\Delta_\mu+W_1 -c(W_2-q)$ is critical on $(M, d\mu)$.
\end{proposition}

\noindent\textbf{Proof:} By Lemma \ref{hpreservescrit}, $\Delta_\nu+W_1-q$ is subcritical on $(M, d\nu)$. Note that $W_1-q, W_2-q\in K^\infty(M, d\nu)$ as well. Since $(M, d\nu)$ satisfies (HKE), by Theorem \ref{zhaomain}, there exists a profile $g>0$ for $\Delta_\nu+W_1-q$ such that $g\asymp 1$.\\

Let $d\theta=g^2d\nu$. Since $(M, d\nu)$ satisfies (HKE), so does $(M, d\theta)$. In fact, by Lemma \ref{greenhke} and the fact that $V_\nu(x,r)\asymp V_\theta(x,r)$ over $x\in M$ and $r>0$, we have that $G_\nu(x,y)\asymp G_\theta(x,y)$. Thus we again have that $W_2-q\in K^\infty(M, d\theta)$. In particular, $W_2-q$ is Green-bounded on $(M, d\theta)$, so there exists $c'>0$ such that $$\sup_{x\in M}\int_{M}G_\theta(x,y)\abs{c'(W_2-q)(y)}d\theta(y)<1.$$ By Lemma \ref{greensmallgauge}, $c'(W_2-q)$ is gaugeable and thus $\Delta_\theta-c'(W_2-q)$ is subcritical. I claim that there exists $c>0$ such that $\Delta_\theta-c(W_2-q)\not\geq 0$. Indeed, by hypothesis there exists $\mc{O}\subseteq M$ open and $\epsilon>0$ such that $W_2>q+\epsilon$ on $\mc{O}$. Let $\phi\in C_c^\infty(M)$ be such that $\phi\geq 0$, $\phi$ is not identically zero and $\textrm{supp}(\phi)\subseteq \mc{O}$. Then picking $c>0$ large such that $$\int_{\mc{O}}(\Delta_\theta\phi)\phi d\theta < c\epsilon\int_{\mc{O}}\phi^2d\theta\leq \int_{\mc{O}} c(W_2-q)\phi^2d\theta$$ implies that $\Delta_\theta-c(W_2-q)\not\geq 0$, proving the claim.\\

Thus we can let $c=\sup\set{c'\colon \Delta_\theta-c'(W_2-q)\geq 0}$. Clearly $\Delta_\theta-c(W_2-q)\geq 0$. If $\Delta_\theta-c(W_2-q)$ is subcritical, then repeating the proof leading up to the claim using the $h$-transform as well as gaugeability implies that $\Delta_\theta-(c+\epsilon)(W_2-q)$ is subcritical for some $\epsilon>0$. But this is impossible. Thus $\Delta_\theta-c(W_2-q)$ is critical.\\

Unraveling, by Lemma \ref{hpreservescrit} we see that $\Delta_\nu + W_1-q-c(W_2-q)$ is critical and hence $\Delta_\mu+W_1-c(W_2-q)$ is critical, as desired.\qed\\

\noindent\textbf{Remark:} If for example we are under the hypotheses of Theorem \ref{quaddecayapplication}, then if we take $W_1, W_2\geq 0$ with $\abs{W_1(x)}, \abs{W_2(x)}\preceq \brac{x}^{-(2+\epsilon)}$ with $W_1$ and $W_2$ not identically zero, then we can take any $q\in C_c^\infty(M)$ such that $q\geq 0$ for which there exists $x\in M$ with $0<q(x)<W_2(x)$ and apply the above to get $c>0$ such that $\Delta_\mu+W_1-c(W_2-q)$ is critical.

\subsection{More Examples}

\begin{example}\label{halfcylinder} The Infinite Half-Cylinder\end{example} For $n\geq 2$ let $S^n=\set{x\in \Real^{n+1}\colon \abs{x}=1}$ be the $n$-dimensional sphere, with Riemannian metric inherited from $\Real^{n+1}$. Let $Q^n=\set{x=(x_1, \dots, x_n)\in S^n\colon x_1\leq 0}$ be the closed hemisphere. Define a manifold $M=Q^{n+1}\cup (0, \infty)\times S^n$. This is not necessarily smooth at the juncture between $Q^{n+1}$ and $(0, \infty)\times S^n$, but this can be smoothed out, yielding a Riemannian manifold. Let $d\mu$ be the canonical Riemannian volume, making $(M, d\mu)$ a weighted manifold.\\

Since $M$ has nonnegative Ricci curvature outside a compact set, $(M, d\mu)$ satisfies (HKE). (See citation.) Let $o\in M$ in the hemisphere be a distinguished point. Note that $M$ satisfies (RCA) with respect to $o$. We have that $V_\mu(r)\asymp r^n$ for $0<r<1$ and $V_\mu(r)\asymp r$ for $r\geq 1$. Also $V_\mu(r)\asymp V_\mu(x,r)$ uniformly in $x\in M$ and $r>0$. Since $\int_{1}^\infty \frac{dt}{V_\mu(\sqrt{t})}=+\infty$, $(M, d\mu)$ is parabolic.\\

In particular we have that \begin{equation}\widehat{H}(y)\asymp H(\abs{y})\asymp 1+\left(\int_1^{\abs{y}}\frac{sds}{V_\mu(s)}\right)_+ \asymp 2+\abs{y}=\brac{y}.\end{equation} Hence if $q\in C_c^\infty(M)$ is nonnegative and not identically zero, and $h>0$ is a profile for $\Delta_\mu+q$, we have \begin{equation}
h(x)\asymp \brac{x}.
\end{equation}
Thus if $W\in C^\infty(M)$ is not identically zero such that $0\leq W(x) \preceq\brac{x}^{-(2+\epsilon)}$ for some $\epsilon>0$, then if $p_\mu^W(t,x,y)$ is the heat kernel of $\Delta_\mu+W$ on $(M, d\mu)$ we get the heat kernel estimate $$\frac{c_1\brac{x}\brac{y}}{\min(t^{n/2}, t^{1/2})(\brac{x}+\sqrt{t})(\brac{y}+\sqrt{t})}\exp\left(-c_2\frac{d(x,y)^2}{t}\right)\leq p_\mu^W(t,x,y)$$ $$\leq\frac{c_3\brac{x}\brac{y}}{\min(t^{n/2}, t^{1/2})(\brac{x}+\sqrt{t})(\brac{y}+\sqrt{t})}\exp\left(-c_4\frac{d(x,y)^2}{t}\right)$$ where $c_1, c_2, c_3, c_4>0$ are constants and $t>0$, $x,y\in M$ are arbitrary.\\

The previous example may be viewed as a model manifold, which we define now. Our presentation is largely the same as that in Grigor'yan and Saloff-Coste \cite{grigoryansaloffcoste2}.

\begin{example}\label{modelmanifolds}
Model Manifolds
\end{example}
A \textbf{model manifold} is $\Real^N$, $N\geq 2$, with a Riemannian metric given in polar coordinates $(r, \theta)\in (0, +\infty)\times S^{N-1}$ by \begin{equation}
ds^2=dr^2+\psi(r)^2d\theta^2
\end{equation}
where $d\theta^2$ is the standard metric on $S^{N-1}$ and $\psi$ is a smooth positive function on $(0, +\infty)$. The metric $ds^2$ can be extended smoothly to all of $\Real^N$ (i.e. across the origin) under the conditions that $\psi(0)=0, \psi'(0)=1$, and $\psi''(0)=0$. Given such a $\psi$, we let $M_\psi$ denote the model manifold $\Real^N$ equipped with the metric $ds^2$. Letting $d\mu$ be the Riemannian volume measure on $M_\psi$, we obtain a weighted manifold $(M_\psi, d\mu)$.\\

Let $o\in M_\psi$ be the origin so that $V_\mu(r)=V_\mu(o, r)$. Due to the radial nature of the metric, $(M_\psi, d\mu)$ satisfies (RCA) with respect to $o$. We have that \begin{equation}\label{modelvol}V_\mu(r)=\omega_N\int_0^r \psi(s)^{N-1}ds,\end{equation} where $\omega_N>0$ is a constant depending on $N$. The volume of a remote ball is given by \begin{equation}V_\mu(x,r)\asymp \begin{cases} 
      r^N, & r\leq \psi(\abs{x}) \\
      r\psi(\abs{x})^{N-1}, & r\geq \psi(\abs{x})
   \end{cases}
\end{equation}

By Proposition 4.10 in Grigor'yan and Saloff-Coste \cite{grigoryansaloffcoste2}, if $$\sup_{[r,2r]}\psi\preceq \inf_{[r,2r]}\psi,$$ $$\psi(r)\preceq r,$$ and $$\int_0^r\psi(s)^{N-1}ds\preceq r\psi(r)^{N-1}$$ over all $r>0$, then $(M_\psi, d\mu)$ satisfies (HKE). In particular, if $\psi(r)\asymp r^\beta$ for $r>1$, then $(M_\psi, d\mu)$ satisfies (HKE) if and only if $-1/(N-1)<\beta\leq 1$.\\

The infinite half-cylinder Example \ref{halfcylinder} is analogous to $\beta=0$. The case $\beta=1$ is analogous to $\Real^N$ with its Euclidean metric. By (HKE), $(M_\psi, d\mu)$ is parabolic if and only if $\int_1^\infty\frac{dt}{V_\mu(\sqrt{t})}=+\infty$, so using Equation \ref{modelvol}, if $\psi(r)\asymp r^\beta$ for $r>1$ with $-1/(N-1)<\beta\leq 1$, then $(M_\psi,d\mu)$ is parabolic if and only if $\beta\leq 1/(N-1)$. In fact we have $V_\mu(r)\asymp r^{\beta(N-1)+1}$ for $r>1$. We assume from here on out that $-1/(N-1)<\beta\leq 1/(N-1)$.\\

For a remote ball $B(y, r)$ with $r>1$ and $\abs{y}>1$, if $r\geq \psi(\abs{y})\asymp \abs{y}^\beta$, then $V_\mu(y, r)\asymp r\abs{y}^{\beta(N-1)}$. If $r\leq \psi(\abs{y})$, then $V_\mu(y, r)\asymp r^N$. \\

Thus we have that for $\abs{y}>1$, $$\widehat{H}(y)=1+\left(\int_1^{\abs{y}^2} \frac{dt}{V_\mu(y, \sqrt{t})}\right)_+$$ $$ \asymp 1+\int_1^{\max(1,\abs{y}^{2\beta})} \frac{dt}{t^{N/2}} + \int_{\max(1,\abs{y}^{2\beta})}^{\abs{y}^2}\frac{dt}{\sqrt{t}\abs{y}^{\beta(N-1)}}\asymp \abs{y}^{1-\beta(N-1)}.$$ Similarly, if $\beta<1/(N-1)$, then $$H(\abs{y})\asymp \brac{y}^{1-\beta(N-1)}\asymp \widehat{H}(y).$$ If $\beta=1/(N-1)$, then $$H(\abs{y})\asymp \log\brac{y}.$$ In both cases we have \begin{equation}\label{modelmainhyp1}\widehat{H}(y)\preceq H(\abs{y})\end{equation} over all $y\in M$. Notably, if $\beta<0$, then for a fixed $r>1$ we have that $$\lim_{\abs{y}\rightarrow \infty} \frac{V_\mu(y, r)}{V_\mu(r)}=0$$ and yet we still have $\widehat{H}(y)\preceq H(\abs{y})$.\\


Due to $\beta\leq 1/(N-1)$, we have that for $1\leq s<r$ we have \begin{equation}\label{modelmainhyp2}\frac{V_\mu(r)}{V_\mu(s)}\preceq \left(\frac{r}{s}\right)^2.\end{equation} 


Putting it all together, we have that: for $-1/(N-1)<\beta\leq 1/(N-1)$ with $\psi(r)\asymp r^\beta$ for $r>1$, the model manifold $(M_\psi, d\mu)$ is complete, non-compact, parabolic, satisfies (HKE) as well as (RCA) with respect to the origin $o$. By \ref{modelmainhyp1} and \ref{modelmainhyp2}, if $W\in C^\infty(M)$ is not identically zero with $0\leq W(x)\preceq \brac{x}^{-(2+\epsilon)}$ for some $\epsilon>0$, then for $\beta=1/(N-1)$, letting $p_\mu^W(t,x, y)$ be the heat kernel of $\Delta_\mu+W$, we get the constants $c_1, c_2, c_3, c_4>0$ such that $$\frac{c_1\log\brac{x}\log \brac{y}}{\log^2(\brac{x}+\sqrt{t})V_\mu(x, \sqrt{t})}e^{-c_2\frac{d(x,y)^2}{t}}\leq p^W(t,x,y)\leq \frac{c_3 \log\brac{x}\log \brac{y}}{\log^2(\brac{x}+\sqrt{t})V_\mu(x, \sqrt{t})}e^{-c_4\frac{d(x,y)^2}{t}}.$$ For $\beta<1/(N-1)$, we instead get $c_1, c_2, c_3, c_4>0$ such that $$\frac{c_1\brac{x}^\sigma\brac{y}^\sigma}{(\brac{x}+\sqrt{t})^{2\sigma} V_\mu(x, \sqrt{t})}e^{-c_2\frac{d(x,y)^2}{t}}\leq p_\mu^W(t,x,y)\leq \frac{c_3\brac{x}^\sigma\brac{y}^\sigma}{(\brac{x}+\sqrt{t})^{2\sigma} V_\mu(x, \sqrt{t})}e^{-c_4\frac{d(x,y)^2}{t}}$$ for all $t>0$ and $x,y\in M_\psi$, where $\sigma=1-\beta(N-1)$. One can write further expand using our bounds on $V_\mu(x,\sqrt{t})$ if desired.\\

\begin{example}\label{almostnonparabolicex}
A Parabolic Manifold with Faster Than Quadratic Volume Growth
\end{example}

Consider the manifold $\Real^2$ with weighted measure $d\mu=\log\brac{x}dx$, where $dx$ is Lebesgue measure. Let $o\in \Real^2$ be the origin so that $(\Real^2, d\mu)$ satisfies (RCA) with respect to $o$. Note first that $V_\mu(r)\asymp r^2\log(2+r)$, so that volume grows faster than quadratically.\\

If $B(x, r)$ is a remote ball then $\log\brac{y}\asymp \log\brac{x}$ for all $y\in B(x, r)$. This immediately yields (VD) and (PI) for remote balls in $(M, d\mu)$. Volume comparison also follows from the fact that if $r=\abs{x}$, then $V_\mu(r)\asymp r^2\log(2+r)=r^2\log \brac{x} \asymp V_\mu(x, r/64)$. Therefore by Theorem \ref{volumecomptheorem}, $(M, d\mu)$ satisfies (HKE).\\

Since $V_\mu(r)\preceq V_\mu(y, r)$ for $0<r\leq \abs{y}$, we immediately have that $\widehat{H}(y)\preceq H(\abs{y})$. Next we have that $$H(\abs{y})=1+\left(\int_1^{\abs{y}^2} \frac{dt}{V_\mu(\sqrt{t})}\right)_+ \asymp 1+\left(\int_{1}^{\abs{y}^2} \frac{dt}{t\log(2+\sqrt{t})}\right)_+\asymp \log\log \brac{y}.$$ Since $H(\abs{y})\rightarrow +\infty$ as $\abs{y}\rightarrow \infty$, $(M, d\mu)$ is parabolic. Note that $V_\mu(x, r)\asymp r^2\log(\brac{x}+r)$. Considering all of the above, if $W\in C^\infty(M)$ is not identically zero with $0\leq W(x)\preceq \brac{x}^{-(2+\epsilon)}$, then we can take $0<\delta<\epsilon$ and note that for $1\leq s<r$, $$\frac{V_\mu(r)}{V_\mu(s)}\preceq \left(\frac{r}{s}\right)^{2+\delta},$$ with constant depending on $\delta$. Therefore we get constants $c_1, c_2, c_3, c_4>0$ such that, letting $p_\mu^W(t,x,y)$ be the heat kernel of $\Delta_\mu+W$, $$\frac{c_1 \log\log\brac{x}\log\log\brac{y}}{t\log(\brac{y}+\sqrt{t})(\log\log(\brac{y}+\sqrt{t}))^2}e^{-c_2\frac{\abs{x-y}^2}{t}}\leq p_\mu^W(t,x,y)\leq$$ $$\frac{c_3 \log\log\brac{x}\log\log\brac{y}}{t\log(\brac{y}+\sqrt{t})(\log\log(\brac{y}+\sqrt{t}))^2}e^{-c_4\frac{\abs{x-y}^2}{t}}$$ for all $t>0$ and $x,y\in \Real^2$.

\subsection{Acknowledgements}

The authors would like to thank Yehuda Pinchover for helpful comments on the manuscript. The first author was supported in part by the National Science Foundation Graduate Research Fellowship grants number DGE-2139899 and DGE-1650441. The second author was supported in part by National Science Foundation grants DMS-2054593 and DMS-2343868.





\bibliographystyle{plain}

\bibliography{ParabolicBib.bib}

\end{document}